\documentclass[leqno,english]{elsarticle}
\usepackage{hyperref}
\usepackage{amssymb}
\usepackage{amsmath}
\usepackage[square,numbers]{natbib}
\usepackage[all]{xy}

\newcommand{\mypart}[1]{\section*{#1}\addcontentsline{toc}{part}{#1}}

\renewcommand{\thesection}{\arabic{section}}
\renewcommand{\thesubsection}{\thesection.\arabic{subsection}}
\renewcommand{\thesubsubsection}{\thesubsection.\arabic{subsubsection}}

\renewcommand{\theenumi}{\arabic{enumi}}
\renewcommand{\labelenumi}{(\theenumi)}

\newtheorem{mainthm}{Theorem}
\newtheorem{mainfact}[mainthm]{Fact}
\newtheorem{mainprop}[mainthm]{Proposition}
\numberwithin{mainthm}{section}

\newtheorem{thm}[subsubsection]{Theorem}
\newtheorem{lemm}[subsubsection]{Lemma}
\newtheorem{prop}[subsubsection]{Proposition}

\newtheorem{fact}[subsubsection]{Fact}
\newtheorem{obsv}[subsubsection]{Observation}

\newdefinition{construction}[subsubsection]{Construction}
\newdefinition{recollection}[subsubsection]{Recollection}
\newdefinition{remark}[subsubsection]{Remark}
\newproof{proof}{Proof}

\DeclareMathOperator{\Mor}{Mor}
\DeclareMathOperator{\End}{\mathrm{End}}
\DeclareMathOperator{\Hom}{\mathrm{Hom}}

\DeclareMathOperator{\Func}{\mathcal{F}}
\DeclareMathOperator{\Op}{\mathcal{O}}

\DeclareMathOperator{\unit}{\mathbf{1}}
\DeclareMathOperator{\EM}{\mathit{EM}}

\DeclareMathOperator{\NN}{\mathbb{N}}
\DeclareMathOperator{\ZZ}{\mathbb{Z}}
\DeclareMathOperator{\kk}{\Bbbk}
\DeclareMathOperator{\FF}{\mathbb{F}}
\DeclareMathOperator{\QQ}{\mathbb{Q}}

\DeclareMathOperator{\AOp}{\mathsf{A}}

\DeclareMathOperator{\COp}{\mathsf{C}}
\DeclareMathOperator{\EOp}{\mathsf{E}}
\DeclareMathOperator{\IOp}{\mathsf{I}}
\DeclareMathOperator{\KOp}{\mathsf{K}}
\DeclareMathOperator{\POp}{\mathsf{P}}
\DeclareMathOperator{\QOp}{\mathsf{Q}}
\DeclareMathOperator{\ROp}{\mathsf{R}}
\DeclareMathOperator{\SOp}{\mathsf{S}}
\DeclareMathOperator{\Free}{\mathsf{F}}

\DeclareMathOperator{\Sym}{\mathit{S}}
\DeclareMathOperator{\Tens}{\mathit{T}}

\DeclareMathOperator{\C}{\mathcal{C}}

\DeclareMathOperator{\E}{\mathcal{E}}
\DeclareMathOperator{\Simp}{\mathcal{S}}
\DeclareMathOperator{\M}{\mathcal{M}}

\DeclareMathOperator{\Id}{Id}
\DeclareMathOperator{\id}{id}

\DeclareMathOperator{\colim}{colim}

\DeclareMathOperator{\Tor}{Tor}

\title{The Bar Complex of an E-infinity Algebra}
\author[bf]{Benoit Fresse}
\fntext[bf]{Research supported in part by ANR grant JCJC06 OBTH}
\ead[bf]{Benoit.Fresse@math.univ-lille1.fr}
\ead[url]{http://math.univ-lille1.fr/\~{ }fresse}
\address{UMR 8524 du CNRS et de l'Universit\'e de Lille 1, Sciences et Technologies\\
Cit\'e Scientifique -- B\^atiment M2\\
F-59655 Villeneuve d'Ascq C\'edex (France)}
\date{August 7, 2009}

\begin{document}

\begin{abstract}
The standard reduced bar complex B(A) of a differential graded algebra A
inherits a natural commutative algebra structure if A is a commutative algebra.
We address an extension of this construction in the context of E-infinity algebras.
We prove that the bar complex of any E-infinity algebra can be equipped with the structure of an E-infinity algebra
so that the bar construction defines a functor from E-infinity algebras to E-infinity algebras.
We prove the homotopy uniqueness of such natural E-infinity structures on the bar construction.

We apply our construction to cochain complexes of topological spaces,
which are instances of E-infinity algebras.
We prove that the n-th iterated bar complexes of the cochain algebra of a space X
is equivalent to the cochain complex of the n-fold iterated loop space of X,
under reasonable connectedness, completeness and finiteness assumptions on X.
\end{abstract}


\maketitle

\mypart{Introduction}

This paper is concerned with the standard reduced bar complex $B(A)$
defined basically for an associative differential graded algebra $A$
equipped with an augmentation over the ground ring $\kk$.
We also consider the natural extension of the bar construction
to $A_\infty$-algebras,
differential graded algebras
equipped with a set of coherent homotopies
that make the structure associative in the strongest homotopical sense.

By a classical construction,
the bar complex of an associative and commutative algebra inherits a multiplicative structure,
unlike the bar complex of a non-commutative algebra,
and still forms a differential graded associative and commutative algebra.
In this paper, we address a generalization of this construction to~$E_\infty$-algebras (E-infinity algebras in plain words),
the notion, parallel to the notion of an $A_\infty$-algebra,
which models a differential graded algebra
equipped with a set of coherent homotopies
that make the structure associative and commutative in the strongest homotopical sense.
Our main theorems, Theorems~\ref{BarStructure:Existence:FunctorLevel}-\ref{BarStructure:Uniqueness:FunctorLevel},
give the existence and the homotopy uniqueness of an $E_\infty$-algebra structure on the bar construction
so that:
\begin{enumerate}
\item
The bar construction $B(A)$
defines a functor from $E_\infty$-algebras to $E_\infty$-algebras.
\item
The $E_\infty$-algebra structure of $B(A)$
reduces to the standard commutative algebra structure of the bar construction
whenever~$A$ is a commutative algebra.
\end{enumerate}

\medskip
To make these assertions more precise,
a model of the category of~$E_\infty$-algebras
has to be fixed.
For this purpose,
we use that the algebra structures which occur in our problem
are modeled by operads:
an~$A_\infty$-algebra is equivalent to an algebra
over an $A_\infty$-operad,
in our context a differential graded operad weakly-equivalent to the operad of associative algebras;
an~$E_\infty$-algebra is an algebra
over an~$E_\infty$-operad,
a differential graded operad weakly-equivalent to the operad of associative and commutative algebras
(see~\cite{May} for the original definition in the topological framework).
Our existence and uniqueness theorems give a functorial $E_\infty$-algebra structure on the bar construction,
for every category of algebras over an $E_\infty$-operad $\EOp$,
for any $E_\infty$-operad $\EOp$.
To define the action on the target,
we just have to take a cofibrant replacement of $\EOp$
with respect to the model structure of differential graded operads~\cite{BergerMoerdijk,Hinich}.

The overall idea of our construction is to use modules over operads to represent functors on categories of algebras over operads.
The bar construction itself is determined by a right module over a particular $A_\infty$-operad,
the chain operad of Stasheff's associahedra (\emph{Stasheff's operad} for short).
The existence and uniqueness of $E_\infty$-algebra structures on the bar construction
is proved at the module level
by techniques of homotopical algebra.
The arguments rely on the existence of a model structure for right modules over operads.

The existence of a dual $E_\infty$-coalgebra structure on the cobar construction
has already been obtained
by a different method in~\cite{JRSmithOperadic}.
But:
the modelling of functors by modules over operads makes our construction more conceptual;
our uniqueness theorem makes the definition of an $E_\infty$-structure easier
since a simple characterization ensures us to obtain the right result.

\medskip
Since the bar construction defines a functor from $E_\infty$-algebras to $E_\infty$-algebras,
we have a well-defined iterated bar complex~$B^n(A)$
associated to any $E_\infty$-algebra.
Our motivation, explained next, is to have an iterated bar complex~$B^n(C^*(X))$,
for any cochain algebra $C^*(X)$, for every pointed topological space $X$,
so that $B^n(C^*(X))$
is equivalent, under reasonable finiteness and connectedness assumptions on the space $X$,
to $C^*(\Omega^n X)$,
the cochain algebra of the iterated loop space $\Omega^n X$.

The usual cochain complexes $C^*(X)$ associated to topological spaces
are examples of objects equipped with an~$E_\infty$-algebra structure
(see~\cite{HinichSchechtman} and the more combinatorial constructions of~\cite{BergerFresse,McClureSmith}).
In positive characteristic,
the existence of Steenrod operations represents a primary obstruction to the existence
of a genuine commutative algebra equivalent to $C^*(X)$
and
one has to use $E_\infty$-algebras (or equivalent notions)
to model faithfully the homotopy of a space $X$
by a cochain complex (see~\cite{Mandell}).

According to classical results of Adams~\cite{AdamsCobar} and Adams-Hilton~\cite{AdamsHilton},
the bar complex $B(C^*(X))$,
where $C^*(X)$ is the cochain algebra of a topological space $X$,
is equivalent as a chain complex to $C^*(\Omega X)$,
the cochain complex of the loop space $\Omega X$.
Since the cochain complex $C^*(X)$ forms an $E_\infty$-algebra,
we obtain by our structure theorem
that the bar complex $B(C^*(X))$ comes equipped with a well-defined $E_\infty$-algebra structure.
To obtain the topological interpretation of the iterated bar complex $B^n(C^*(X))$,
we prove that $B(C^*(X))$ is equivalent to $C^*(\Omega X)$
as an $E_\infty$-algebra.

For this aim,
we prove that, for a cofibrant $E_\infty$-algebra,
the usual bar construction
is equivalent as an $E_\infty$-algebra to a categorical version of the bar construction
in which tensor products are replaced by algebra coproducts.
Then
we apply a theorem of Mandell~\cite{Mandell} which asserts that the categorical bar construction
of a cofibrant replacement of $C^*(X)$
defines an $E_\infty$-algebra equivalent to $C^*(\Omega X)$.

The categorical bar construction preserves weak-equivalences between cofibrant $E_\infty$-algebras only.
Therefore we have to form a cofibrant replacement of $C^*(X)$ in $E_\infty$-algebras
in order to apply the categorical bar construction reasonably.
In contrast,
the usual bar construction preserves weak-equivalences between all $E_\infty$-algebras
which are cofibrant in the underlying category of dg-modules (all $E_\infty$-algebras if the ground ring is a field).
For this reason,
we can apply the usual bar construction
to the cochain algebra itself $C^*(X)$, and not only to a cofibrant replacement of $C^*(X)$,
to still have an $E_\infty$-algebra equivalent to $C^*(\Omega X)$.

The article~\cite{Mandell}
gives an attractive theoretical setting to model the homotopy of spaces in positive characteristic,
but in practice one has to face deep difficulties to build cofibrant replacements in categories of $E_\infty$-algebras.
In this sense,
our construction gives an effective substitute for the categorical bar construction
used in~\cite{Mandell}.

According to~\cite{Miller},
the bar complex of simplicial commutative algebras
models the suspension in the homotopy category of simplicial commutative algebras.
In passing,
we prove that, in the differential graded setting,
the bar complex of $E_\infty$-algebras yields a model of the suspension
in categories of $E_\infty$-algebras.

\medskip
Other attempts to define an iterated bar construction
occur in the literature
outside J.R.Smith's memoirs~\cite{JRSmith,JRSmithOperadic}.
Usually,
authors deal with the dual cobar construction and chain complexes rather than cochain complexes.
If we assume reasonable finiteness assumptions on spaces,
then this dual construction is equivalent to the bar construction
and nothing changes.
To simplify we examine the previous results of the literature
in the context of the bar construction.
\begin{enumerate}
\item\label{GeometricCobar}
The original geometrical approach of Adams~\cite{AdamsCobar} and Adams-Hilton~\cite{AdamsHilton}
is continued by Milgram in~\cite{Milgram} and Baues in~\cite{BauesDoubleBar,BauesBarGeometry,BauesCobarHopf}
to define a double bar construction $B^2(C^*(X))$,
for any cochain algebra $C^*(X)$, where $X$ is a simplicial set (see also the survey article~\cite{CarlssonMilgram}).
\item\label{PerturbativeCobar}
In~\cite{KadeishviliSaneblidze},
Kadeishvili-Saneblidze use perturbation lemmas and the classical chain equivalence $B(C^*(X))\sim C^*(\Omega X)$
to obtain an inductive construction of an iterated bar complex $B^n(C^*(X))$
together with a chain equivalence $B^n(C^*(X))\sim C^*(\Omega^n X)$,
for every cochain algebra $C^*(X)$;
this approach is used by Rubio-Sergeraert in the Kenzo program~\cite{RubioSergeraert}
to perform computer calculations.
\item\label{NonCommutativeGeometricCobar}
In~\cite{Karoubi},
Karoubi uses ideas of non-commutative differential geometry
and non-commutative analogues of difference calculus
to introduce new cochain complexes~$D^*(X)$
for which a modified iterated bar complex $B^n(D^*(X))$ can be defined
so that $B^n(D^*(X))\sim D^*(\Omega^n X)$.
\end{enumerate}

The difficulty in~(\ref{GeometricCobar}) is to understand the geometry of certain cell complexes
in order to define higher iterated bar complex $B^n(C^*(X))$ for $n>2$ (see~\cite{BauesBarGeometry,BilleraKapranovSturmfels}).
In the approach of~(\ref{GeometricCobar}), and similarly in~(\ref{NonCommutativeGeometricCobar}),
the bar construction is only defined for complexes of a particular type.
In the approach of~(\ref{PerturbativeCobar}), one has to keep track of a simplicial model of $\Omega^n X$,
the iterated Kan construction~$G^n(X)$,
to define the differential of $B^n(C^*(X))$.

In contrast,
our theorems imply the existence of a well-characterized iterated bar complex $B^n(A)$,
for every $E_\infty$-algebra $A$,
and such that $B^n(A)$ incorporates minimal information in itself.
Besides,
we have to use multiplicative structures to relate the iterated bar complex $B^n(C^*(X))$
to the cochain complex of an iterated loop space $C^*(\Omega^n X)$,
but the iterated bar complex $B^n(A)$ can be determined directly by using that a composite of functors associated to modules over operads,
like the iterated bar complex,
forms itself a functor determined by a module over an operad (see~\cite{FresseModules}).
This observation,
beyond the scope of this article,
is the starting point of~\cite{FresseIteratedBar}.

\medskip
In this introduction,
we adopt the usual convention to apply the bar construction to augmented unital algebras.
In the context of the cochain complex of a space $X$,
the augmentation is determined by the choice of a base point $*\in X$.
But in the definition of the bar complex we have to replace an algebra $A$ by its augmentation ideal $\bar{A}$,
which forms a non-augmented non-unital algebra,
and the cochain complex of a space $C^*(X)$
by the associated reduced complex $\bar{C}^*(X)$.
Therefore
it is more natural to use non-augmented non-unital algebras for our purpose
and we take this convention in the core sections of the article
(for details, see~\S\ref{BarConstruction:AinfinityBarComplex:Definition} and \S\ref{BarConstruction:AinfinityBarComplex:UnitaryContext}).

\mypart{Contents}

In the first part of the paper, \emph{``Background''},
we survey new ideas introduced in~\cite{FresseModules}
to model functors on algebras over operads
by modules over operads.
These preliminaries are necessary to make the conceptual setting of our constructions
accessible to readers which are only familiar with standard definitions
of the theory of operads.

The object of our study, the bar construction,
appears in the second part, \emph{``The bar construction and its multiplicative structure''},
where we prove the main results of the article.
In the core sections, \S\S\ref{BarConstruction}-\ref{HomotopyInterpretation},
we define the bar module,
the module over Stasheff's operad which represents the bar construction,
we prove the existence and uniqueness of a multiplicative structure on the bar construction,
and we give a homotopy interpretation of the bar construction
in the model category of $E_\infty$-algebras.
For a more detailed outline,
we refer to the introduction of this part.

In the concluding part, \emph{``The iterated bar construction and iterated loop spaces''},
we address topological applications of our results.
As explained in this introduction,
we use the multiplicative structure of the bar construction
to define an iterated bar construction $B^n(C^*(X))$,
for any cochain algebra $C^*(X)$,
so that $B^n(C^*(X))\sim C^*(\Omega^n X)$.
One aim of this part is to make explicit reasonable finiteness, completeness and connectedness assumptions on $X$
which ensure this equivalence.

\renewcommand{\themainthm}{\thesection.\Alph{mainthm}}

\mypart{Background}

Before studying the structure of the bar construction,
we survey ideas introduced in the book~\cite{FresseModules} to make the overall setting of our constructions accessible to readers.

First,
our use of functors and modules over operads motivates a review of the categorical background of operad theory,
to which~\S\ref{Background:SymmetricMonoidalCategories}
are devoted.
Then,
in~\S\ref{Background:OperadAlgebras},
we review the definition of an operad, of an algebra over an operad,
and the definition of categories of modules associated to operads.
The correspondence between modules over operads and functors is addressed in~\S\ref{Background:FunctorModules}.

Throughout the paper,
we use extensively extension and restriction functors in the context of algebras and modules over operads.
The last subsection of this part,~\S\ref{Background:ExtensionRestriction},
is devoted to recollections
on these topics.

\subsection{Symmetric monoidal categories over dg-modules}\label{Background:SymmetricMonoidalCategories}
As usual in the literature,
we assume that operads consist of objects in a fixed base symmetric monoidal category
-- for our purpose,
the category of unbounded differential graded modules (\emph{dg-modules} for short)
over a fixed ground ring $\kk$ (see~\S\ref{Background:SymmetricMonoidalCategories:Definition}).

In contrast,
we can assume that the underlying category of algebras over an operad is not the base category itself, to which the operad belongs,
but some symmetric monoidal category over the category of dg-modules.
Though we only use specific examples of such categories in applications, the category of dg-modules itself,
the category of $\Sigma_*$-modules, and categories of right modules over an operad,
for which alternative point of views are available (see~\S\ref{Background:OperadAlgebras}),
we prefer to review the definition of this general setting which gives the right conceptual background
to understand our arguments.

\subsubsection{Symmetric monoidal categories over dg-modules}\label{Background:SymmetricMonoidalCategories:Definition}
Let $\kk$ be a ground ring, fixed once and for all.
Throughout the paper,
the notation $\C$ refers to the category of dg-modules,
where a dg-module consists of a lower $\ZZ$-graded $\kk$-module $C = \oplus_{*\in\ZZ} C_*$
equipped with an internal differential, usually denoted by $\delta: C\rightarrow C$,
that decreases degrees by $1$.
The usual convention $C^* = C_{-*}$ makes any upper graded module equivalent to an object of $\C$.

The category of dg-modules is equipped with the standard tensor product of dg-modules $\otimes: \C\times\C\rightarrow\C$
which provides $\C$ with the structure of a symmetric monoidal category.
The unit object of dg-modules is formed by the ground ring itself $\kk$,
viewed as a dg-module concentrated in degree $0$.

For us,
a symmetric monoidal category over $\C$ is a symmetric monoidal category~$\E$
equipped with an external tensor product $\otimes: \C\times\E\rightarrow\E$
so that an obvious generalization of relations of symmetric monoidal categories
holds in $\E$,
for any composite of the tensor products $\otimes: \E\times\E\rightarrow\E$ and $\otimes: \C\times\E\rightarrow\E$.
For details on this background we refer to~\cite[\S 1.1]{FresseModules}.

In principle,
we assume that the internal tensor product of $\E$, as well as the external tensor product over dg-modules $\otimes: \C\times\E\rightarrow\E$,
preserves colimits.
Under mild set-theoretic assumptions,
these conditions are equivalent to the existence of right adjoints for the internal tensor product
and the external tensor product of $\E$.
In the paper,
we only use the existence of the external-hom
\begin{equation*}
\Hom_{\E}(-,-): \E^{op}\times\E\rightarrow\C,
\end{equation*}
which satisfies
\begin{equation*}
\Mor_{\E}(C\otimes E,F) = \Mor_{\C}(C,\Hom_{\E}(E,F)),
\end{equation*}
for $C\in\C$, $E,F\in\E$.

\subsubsection{Symmetric monoidal model categories over dg-modules}\label{Background:SymmetricMonoidalCategories:ModelCategories}
The category of dg-modules $\C$
is equipped with a cofibrantly generated model structure
such that a morphism $f: C\rightarrow D$ is a weak-equivalence if $f$ induces an isomorphism in homology,
a fibration if $f$ is degreewise surjective,
and a cofibration if $f$ has the left lifting properties with respect to acyclic fibrations.

This model structure is symmetric monoidal (see~\cite[\S 4]{Hovey})
in the sense that:
\begin{enumerate}\renewcommand{\labelenumi}{MM\arabic{enumi}.}\setcounter{enumi}{-1}
\item
\emph{The unit of the tensor product forms a cofibrant object in $\C$}.
\item
\emph{The tensor product $\otimes: \C\times\C\rightarrow\C$ satisfies the pushout-product axiom --
explicitly:
the natural morphism
\begin{equation*}
(i_*,j_*): A\otimes D\bigoplus_{A\otimes C} B\otimes C\rightarrow B\otimes D
\end{equation*}
induced by cofibrations $i: A\rightarrowtail B$ and $j: C\rightarrowtail D$
forms a cofibration in $\C$,
an acyclic cofibration if $i$ or $j$ is also acyclic}.
\end{enumerate}

In the paper,
we use cofibrantly generated model categories $\E$
which are symmetric monoidal over the base category of dg-modules $\C$
and such that the analogues of axioms MM0-MM1 are satisfied at the level of $\E$:
the unit object $\unit\in\E$ forms a cofibrant object in $\E$
and
the internal tensor product of $\E$, as well as the external tensor product of $\E$ over the category of dg-modules $\C$,
satisfies the pushout product axiom.
In this context,
we say that $\E$ forms a cofibrantly generated symmetric monoidal model category
over dg-modules.

The books~\cite{Hirschhorn,Hovey} are our references on the background of model categories.
For the definition of a symmetric monoidal model category,
we refer more particularly to~\cite[\S 4]{Hovey}.
For the generalization of this notion to our relative setting,
we refer to~\cite[\S 11.3]{FresseModules}.

\subsubsection{Enriched model category structures}\label{Background:SymmetricMonoidalCategories:EnrichedModelCategories}
The axioms of symmetric monoidal model categories
are used implicitly when we define the model category of operads
and the model category of algebras over an operad.
In the article,
we also use a dual version of the pushout-product axiom
which holds for the external hom functor of a symmetric monoidal category over $\C$:
\begin{enumerate}\renewcommand{\labelenumi}{MM\arabic{enumi}'.}
\item
\emph{The natural morphism
\begin{equation*}
\qquad\qquad\Hom_{\E}(B,C)\xrightarrow{(i^*,p_*)}\Hom_{\E}(A,C)\times_{\Hom_{\E}(A,D)}\Hom_{\E}(B,D)
\end{equation*}
induced by a cofibration $i: A\rightarrowtail B$ and a fibration $p: C\twoheadrightarrow D$ forms a fibration in $\C$,
an acyclic fibration if $i$ or $p$ is also acyclic}.
\end{enumerate}
The characterization of (acyclic) fibrations in a model category by the left lifting property with respect to (acyclic) cofibrations
and the definition of the external hom
imply readily that axiom MM1' is formally equivalent to the pushout product axiom MM1
for the external tensor product $\otimes: \C\times\E\rightarrow\E$.

\subsection{Operads, algebras and modules over operads}\label{Background:OperadAlgebras}
In this subsection,
we review basic definitions of the theory of operads in the context of symmetric monoidal categories over dg-modules.
To begin with,
we recall briefly the definition of a $\Sigma_*$-module, of an operad,
and of module structures associated to operads.
For details,
we refer to relevant sections of~\cite{FresseModules}.

\subsubsection{Operads and modules over operads}\label{Background:OperadAlgebras:OperadModules}
Throughout the paper,
we use the notation $\M$ to refer to the category of $\Sigma_*$-objects in dg-modules (\emph{$\Sigma_*$-modules} for short),
whose objects are collections $M = \{M(n)\}_{n\in\NN}$,
where $M(n)$ is a dg-module equipped with an action of the symmetric group in $n$ letters $\Sigma_n$, for $n\in\NN$.

In the classical theory,
a module of symmetric tensors
\begin{equation*}
\Sym(M,E) = \bigoplus_{n=0}^{\infty} (M(n)\otimes E^{\otimes n})_{\Sigma_n}
\end{equation*}
is associated to any $\Sigma_*$-module $M\in\M$.
The coinvariants $(M(n)\otimes E^{\otimes n})_{\Sigma_n}$
identify the natural action of permutations on $E^{\otimes n}$
with their action on~$M(n)$.
For our purpose,
we note that this construction makes sense in any symmetric monoidal category $\E$
over the category of dg-modules $\C$,
so that the map $\Sym(M): E\mapsto\Sym(M,E)$
defines a functor $\Sym(M): \E\rightarrow\E$.

The category of $\Sigma_*$-modules comes equipped with a composition product $\circ: \M\times\M\rightarrow\M$
such that $\Sym(M\circ N,E) = \Sym(M,\Sym(N,E))$,
for all $M,N\in\M$, $E\in\E$,
and for every symmetric monoidal category over dg-modules $\E$.
The composition product of $\Sigma_*$-modules
is associative and unital.
The composition unit is defined by the $\Sigma_*$-module
\begin{equation*}
\IOp(n) = \begin{cases} \kk, & \text{if $n=1$}, \\ 0, & \text{otherwise}, \end{cases}
\end{equation*}
and we have $\Sym(\IOp) = \Id$, the identity functor on $\E$.

There are several equivalent definitions for the notion of an operad.
According to one of them,
an operad consists of a $\Sigma_*$-module $\POp$
equipped with an associative product $\mu: \POp\circ\POp\rightarrow\POp$, the composition product of $\POp$,
together with a unit represented by a morphism $\eta: \IOp\rightarrow\POp$.

The structure of a right module over an operad $\ROp$
is defined by a $\Sigma_*$-module $M$
equipped with a right $\ROp$-action determined by a morphism $\rho: M\circ\ROp\rightarrow M$
which is associative with respect to the operad composition product
and unital with respect to the operad unit.
The category of right $\ROp$-modules is denoted by $\M_{\ROp}$.

There is a symmetrically defined notion of left module over an operad $\POp$
consisting of a $\Sigma_*$-module $N$
equipped with a left $\POp$-action determined by a morphism $\lambda: \POp\circ N\rightarrow N$.
One can also define the notion of a bimodule as a $\Sigma_*$-object $N$
equipped with both a right $\ROp$-action $\rho: N\circ\ROp\rightarrow N$
and a left $\POp$-action $\lambda: \POp\circ N\rightarrow N$
that commute to each other.
The notation ${}_{\POp}\M$ refers to the category of left $\POp$-modules
and the notation ${}_{\POp}\M{}_{\ROp}$ to the category of $\POp$-$\ROp$-bimodules.

Note that an operad $\ROp$ forms obviously a right module (respectively, left module, bimodule)
over itself.

The composition product of $\Sigma_*$-modules is not symmetric since this operation is supposed to represent
the composition of functors.
For this reason,
left and right operad actions on $\Sigma_*$-modules
have different nature
though definitions are symmetrical.
In~\S\ref{Background:OperadAlgebras:GeneralizedAlgebras},
we observe that left modules (respectively, bimodules) over operads are equivalent to algebras over operads
and this equivalent definition reflects the structure of left modules
and bimodules more properly.

\subsubsection{The symmetric monoidal category of $\Sigma_*$-modules}\label{Background:OperadAlgebras:SigmaModules}
The category of $\Sigma_*$-modules, which defines an underlying category for operads and modules over operads,
gives our primary example of a symmetric monoidal model category over dg-modules
(outside the category of dg-modules itself).

The unit of the tensor structure of $\Sigma_*$-modules
is the $\Sigma_*$-module $\unit$
such that $\unit(0) = \kk$ and $\unit(n) = 0$ for $n>0$.
The tensor product $C\otimes N\in\M$ of a $\Sigma_*$-module $M\in\M$
with a dg-module $C\in\C$
is given by the obvious formula
\begin{equation*}
(C\otimes M)(r) = C\otimes M(r),
\end{equation*}
for $r\in\NN$.
The tensor product $M\otimes N\in\M$ of $\Sigma_*$-modules $M,N\in\M$
is defined by a formula of the form:
\begin{equation*}
(M\otimes N)(r) = \bigoplus_{s+t=r} \Sigma_r\otimes_{\Sigma_s\times\Sigma_t} M(s)\otimes N(t),
\end{equation*}
for $r\in\NN$.
At the functor level,
the tensor operations of $\Sigma_*$-modules
represent the pointwise tensor products
\begin{equation*}
\Sym(M\otimes N,E) = \Sym(M,E)\otimes\Sym(N,E)\quad\text{and}\quad\Sym(C\otimes M,E) = C\otimes\Sym(M,E),
\end{equation*}
for any symmetric monoidal category over dg-modules $\E$,
where $E\in\E$ (we refer to~\cite[\S 2.1]{FresseModules} for details on these recollections).

Since $\M$ forms a symmetric monoidal category over dg-modules,
a $\Sigma_*$-module~$M$ gives rise to a functor $\Sym(M): \M\rightarrow\M$
on the category of $\Sigma_*$-modules itself.
In fact,
we have an identity $\Sym(M,N) = M\circ N$,
for all $M,N\in\M$
(see~\cite[\S 2.2]{FresseModules}).

\subsubsection{The symmetric monoidal category of right $\ROp$-modules}\label{Background:OperadAlgebras:RightModules}
According to~\cite[\S 6.1]{FresseModules},
the tensor product $M\otimes N$ of right modules over an operad $\ROp$
inherits the structure of a right $\ROp$-module
and similarly for the external tensor product $C\otimes M$
of a dg-module $C\in\C$ with a right $\ROp$-module $M\in\M_{\ROp}$.
Hence
the category of right $\ROp$-modules forms a symmetric monoidal category over dg-modules
so that the forgetful functor $U: \M_{\ROp}\rightarrow\M$
preserves symmetric monoidal structures.
In the context of right $\ROp$-modules,
the functor $\Sym(M): \M_{\ROp}\rightarrow\M_{\ROp}$ is still given by the formula $\Sym(M,N) = M\circ N$,
for all $M\in\M$, $N\in\M_{\ROp}$,
where $M\circ N$ has an obvious right $\ROp$-action induced by the right $\ROp$-action on $N$.

\subsubsection{Symmetric monoidal model structures}\label{Background:OperadAlgebras:UnderlyingModelCategories}
The category of $\Sigma_*$-modules $\M$ inherits a natural model structure
such that a morphism $f: M\rightarrow N$ is a weak-equivalence (respectively, a fibration)
if the underlying collection of dg-module morphisms $f: M(n)\rightarrow N(n)$
consists of weak-equivalences (respectively, fibrations)
in the category of dg-modules.
Cofibrations are determined by the right lifting property with respect to acyclic fibrations.
The model category $\M$ is also cofibrantly generated
and symmetric monoidal over dg-modules
in the sense of~\S\ref{Background:SymmetricMonoidalCategories:ModelCategories}
(see~\cite[\S 11.4]{FresseModules}).

In~\cite[\S 14]{FresseModules}
we check that the category of right modules over an operad $\ROp$
forms a cofibrantly generated symmetric monoidal model category over dg-modules,
like the category of $\Sigma_*$-modules,
provided that the underlying collection of the operad $\{\ROp(n)\}_{n\in\NN}$
consists of cofibrant dg-modules.
Throughout the paper,
we assume tacitely that an operad $\ROp$ satisfies this condition if we deal with model structures
of the category of right $\ROp$-modules.
As usual,
we assume that a morphism of right $\ROp$-modules $f: M\rightarrow N$ is a weak-equivalence (respectively, a fibration)
if the underlying collection
consists of weak-equivalences (respectively, fibrations) of dg-modules $f: M(n)\rightarrow N(n)$
and we characterize cofibrations by the left lifting property with respect to acyclic fibrations.

\subsubsection{On algebras over operads}\label{Background:OperadAlgebras:GeneralizedAlgebras}
In standard definitions,
one uses that the functor $\Sym(\POp)$ associated to an operad $\POp$ forms a monad
in order to define the category of algebras associated to $\POp$.
The usual definition can readily be extended in the context of symmetric monoidal categories over dg-modules
since according to the construction of~\S\ref{Background:OperadAlgebras:OperadModules}
we have a functor $\Sym(\POp): \E\rightarrow\E$
for every symmetric monoidal category $\E$ over the category of dg-modules $\C$.

The structure of a $\POp$-algebra in $\E$
consists of an object $A\in\E$
equipped with an evaluation morphism $\lambda: \Sym(\POp,A)\rightarrow A$
that satisfies natural associativity and unit relations.
The definition of $\Sym(\POp,A)$
implies that the evaluation morphism is also equivalent
to a collection of equivariant morphisms
\begin{equation*}
\lambda: \POp(n)\otimes A^{\otimes n}\rightarrow A
\end{equation*}
formed in the category $\E$.
Throughout the paper,
we use the notation ${}_{\POp}\E$
to refer to the category of $\POp$-algebras in $\E$.

For $E\in\E$,
the object $\Sym(\POp,E)\in\E$ is equipped with a natural $\POp$-algebra structure
and represents the free object associated to $E$
in the category of $\POp$-algebras.
In the paper,
we use the notation $\POp(E) = \Sym(\POp,E)$ to refer to the object $\Sym(\POp,E)$ equipped with the free $\POp$-algebra structure
and we keep the notation $\Sym(\POp,E)$ to refer to the underlying object in $\E$.

In the case $\E = \M$,
we have an identity $\Sym(\POp,M) = \POp\circ M$
from which we deduce that a $\POp$-algebra in $\Sigma_*$-objects is equivalent to a left $\POp$-module.
In the case $\E = \M_{\ROp}$,
we obtain that a $\POp$-algebra in right $\ROp$-modules is equivalent to a $\POp$-$\ROp$-bimodule.
Our conventions for categories
of algebras over operads
is coherent with the notation of~\S\ref{Background:OperadAlgebras:GeneralizedAlgebras} for the category of left $\POp$-modules~${}_{\POp}\M$
and for the category of $\POp$-$\ROp$-bimodules~${}_{\POp}\M{}_{\ROp}$.

In the paper,
we use repeatedly the observation, made in~\S\ref{Background:OperadAlgebras:OperadModules}, that an operad forms a bimodule over itself,
and hence an algebra over itself in the category of right modules over itself.

The categories of right modules over an operad carry the same structures as usual categories of modules over algebras.
The categories of left modules over an operad, as well as the categories of bimodules,
have
structures of different nature
that the notion of an algebra in a symmetric monoidal category over dg-modules
reflects.
The idea of an algebra in a symmetric monoidal category over dg-modules
is also more natural in constructions of this article.
Therefore, in this paper, we prefer to use the language of algebras in symmetric monoidal categories
for left modules and bimodules over operads.

\subsubsection{Model categories of algebras over operads}\label{Background:OperadAlgebras:AlgebraModelCategories}
Let $\POp$ be a $\Sigma_*$-cofibrant operad,
an operad which forms a cofibrant object in the underlying category of $\Sigma_*$-modules.

Let $\E$ be a cofibrantly generated symmetric monoidal model category over dg-modules.
The category of $\POp$-algebras in $\E$
inherits a semi-model structure
such that a morphism $f: A\rightarrow B$ defines a weak-equivalence (respectively, a fibration) in ${}_{\POp}\E$
if $f$ forms a weak-equivalences (respectively, fibrations)
in the underlying category $\E$
(see~\cite{Spitzweck}, we also refer to~\cite{HoveySemiModel} for the notion of a semi-model category).
Roughly,
all axioms of a model category are satisfied in ${}_{\POp}\E$, including M4 and M5,
as long as the source of the morphism $f: A\rightarrow B$
that occurs in these properties is assumed to be cofibrant.

This assertion can be applied to the category of $\Sigma_*$-modules $\E = \M$
(respectively, to the category of right modules over an operad $\E = \M_{\ROp}$)
to obtain that the left $\POp$-modules ${}_{\POp}\M$
(respectively, the $\POp$-$\ROp$-bimodules ${}_{\POp}\M{}_{\ROp}$)
form a semi-model category.

\subsubsection{Model categories of operads}\label{Background:OperadAlgebras:OperadModelCategory}
The category of operads $\Op$ carries a semi-model structure such that the forgetful functor $U: \Op\rightarrow\M$
creates fibrations and weak-equivalences (see~\cite{Spitzweck}).
Thus,
according to definitions for $\M$,
a morphism $f: \POp\rightarrow\QOp$ is a weak-equivalence (respectively, a fibration) in $\Op$
if the underlying morphisms of dg-modules $f: \POp(n)\rightarrow\QOp(n)$, $n\in\NN$, are all weak-equivalences (respectively, fibration)
in the category of dg-modules.
In the core sections of the paper,
we use operads $\POp$ such that $\POp(0) = 0$.
According to~\cite{BergerMoerdijk,Hinich},
the subcategory $\Op_0\subset\Op$ formed by these operads inherits a full model category structure.

As usual,
we characterize cofibrations by the right lifting property with respect to acyclic fibrations in $\Op$.
In particular,
an operad $\POp\in\Op$
is cofibrant as an operad if the lifting exists
in all diagrams of the form
\begin{equation*}
\xymatrix{ & \ROp\ar@{->>}[d]_{\sim}^{p} \\ \POp\ar@{.>}[ur]^{\exists?}\ar[r] & \SOp },
\end{equation*}
where $p: \ROp\rightarrow\SOp$ is an acyclic fibration of operads.

Recall that an operad $\POp$ is said to be $\Sigma_*$-cofibrant
if $\POp$ forms a cofibrant object in the underlying category of $\Sigma_*$-modules.
One can check that cofibrant operads are $\Sigma_*$-cofibrant (see~\cite[Proposition 4.3]{BergerMoerdijk}),
but the converse assertion does not hold.

\subsection{Modules over operads and functors}\label{Background:FunctorModules}
In this subsection,
we recall the definition and categorical properties of functors associated to right modules over operads.

\subsubsection{The functor associated to a right module over an operad}\label{Background:FunctorModules:Definition}
Let $M$ be a right module over an operad $\ROp$.

Let $\E$ be any symmetric monoidal category over dg-modules.
For an $\ROp$-algebra $A\in{}_{\ROp}\E$,
we form the coequalizer:
\begin{equation*}
\xymatrix{\Sym(M\circ\ROp,A)\ar@<+1mm>[r]^{d_0}\ar@<-1mm>[r]_{d_1} & \Sym(M,A)\ar[r] & \Sym_{\ROp}(M,A) },
\end{equation*}
where $d_0$ is the morphism
\begin{equation*}
\Sym(M\circ\ROp,A)\xrightarrow{\Sym(\rho,A)}\Sym(M,A)
\end{equation*}
induced by the right $\ROp$-action on $M$
and $d_1$ is the morphism
\begin{equation*}
\Sym(M\circ\ROp,A) = \Sym(M,\Sym(\ROp,A))\xrightarrow{\Sym(M,\lambda)}\Sym(M,A)
\end{equation*}
induced by the left $\ROp$-action on $A$.

The map $\Sym_{\ROp}(M): A\mapsto\Sym_{\ROp}(M,A)$
defines the functor $\Sym_{\ROp}(M): {}_{\ROp}\E\rightarrow\E$
associated to $M$.
Let $\Func{}_{\ROp}$ denote the category of functors $F: {}_{\ROp}\E\rightarrow\E$.
The definition of $\Sym_{\ROp}(M)$ is obviously natural in $M$
so that $\Sym_{\ROp}: M\mapsto\Sym_{\ROp}(M)$
defines a functor $\Sym_{\ROp}: \M{}_{\ROp}\rightarrow\Func{}_{\ROp}$.

The definition of the functor $\Sym_{\ROp}(M): {}_{\ROp}\E\rightarrow\E$
can be applied to the category of $\Sigma_*$-modules $\E = \M$,
or to another category of right modules $\E = \M_{\SOp}$,
for any operad $\SOp$.
In this context,
the object $\Sym_{\ROp}(M,N)$ is identified with the classical relative composition product $M\circ_{\ROp} N$
of the operad literature.
Indeed,
the relative composition product $M\circ_{\ROp} N$
is defined by a coequalizer of the same form where the objects $\Sym(M,N)$ are replaced by the equivalent composites $\Sym(M,N) = M\circ N$
in the category of $\Sigma_*$-modules (see for instance~\cite[\S 2.1.7]{FressePartitions} for this definition).

\subsubsection{Categorical operations on functors associated to right-modules over operads}\label{Background:FunctorModules:Operations}
To determine the functor $\Sym_{\ROp}(M): {}_{\ROp}\E\rightarrow\E$
associated to a right $\ROp$-module~$M$,
we essentially use:
\begin{enumerate}
\item\label{Background:FunctorModules:Operations:Constant}
For the unit object $M = \unit$,
the functor $\Sym_{\ROp}(\unit): {}_{\ROp}\E\rightarrow\E$ is the constant functor $\Sym_{\ROp}(\unit,A)\equiv\kk$.
\item\label{Background:FunctorModules:Operations:TensorProduct}
We have a natural isomorphism $\Sym_{\ROp}(M\otimes N,A) = \Sym_{\ROp}(M,A)\otimes\Sym_{\ROp}(N,A)$,
for all $M,N\in\M{}_{\ROp}$, $A\in{}_{\ROp}\E$,
and a natural isomorphism $\Sym_{\ROp}(C\otimes M,A) = C\otimes\Sym_{\ROp}(M,A)$,
for all $C\in\C$, $M\in\M{}_{\ROp}$, $A\in{}_{\ROp}\E$,
so that the map $\Sym_{\ROp}: M\mapsto\Sym_{\ROp}(M)$ defines a functor of symmetric monoidal categories over dg-modules
$\Sym_{\ROp}: (\M{}_{\ROp},\otimes,\unit)\rightarrow(\Func{}_{\ROp},\otimes,\kk)$,
the tensor structure of functors being defined pointwise.
\item\label{Background:FunctorModules:Operations:Colimits}
The functor $\Sym_{\ROp}: \M{}_{\ROp}\rightarrow\Func{}_{\ROp}$ preserves colimits.
\end{enumerate}
We refer to~\cite[\S\S 5-6]{FresseModules}
for the proof of these assertions.

The functor $\Sym_{\ROp}: M\mapsto\Sym_{\ROp}(M)$
is uniquely characterized by~(\ref{Background:FunctorModules:Operations:Constant}-\ref{Background:FunctorModules:Operations:Colimits})
and assertion (\ref{Background:FunctorModules:FunctorsToAlgebras:Identity}) of~\S\ref{Background:FunctorModules:FunctorsToAlgebras}
(use the form of generating objects in $\M_{\ROp}$, see~\cite[\S 7.1]{FresseModules}).
If we forget algebra structures on the target,
then this latter assertion implies:
\begin{enumerate}\setcounter{enumi}{3}
\item\label{Background:FunctorModules:Operations:Forgetful}
The functor $\Sym_{\ROp}(\ROp): {}_{\ROp}\E\rightarrow\E$
associated to the operad $\ROp$, viewed as a right module over itself,
represents the forgetful functor $U: {}_{\ROp}\E\rightarrow\E$.
\end{enumerate}

\subsubsection{On algebras in right-modules over operads and functors}\label{Background:FunctorModules:FunctorsToAlgebras}
The assertions of~\S\ref{Background:FunctorModules:Operations}
imply that the evaluation morphism of a $\POp$-algebra in right $\ROp$-modules
\begin{equation*}
\lambda: \POp(n)\otimes N^{\otimes n}\rightarrow N
\end{equation*}
give rise to natural evaluation morphisms at the functor level:
\begin{equation*}
\POp(n)\otimes\Sym_{\ROp}(N,A)^{\otimes n} = \Sym_{\ROp}(\POp(n)\otimes N^{\otimes n},A)\rightarrow\Sym_{\ROp}(N,A),
\end{equation*}
where $A\in{}_{\ROp}\E$.
Thus we obtain that the map $\Sym_{\ROp}(N): A\mapsto \Sym_{\ROp}(N,A)$
defines a functor $\Sym_{\ROp}(N): {}_{\ROp}\E\rightarrow{}_{\POp}\E$.

According to~\cite[Observation 9.2.2]{FresseModules}:
\begin{enumerate}
\item\label{Background:FunctorModules:FunctorsToAlgebras:Identity}
The identity functor $\Id: {}_{\ROp}\E\rightarrow{}_{\ROp}\E$
is realized by the functor $\Sym_{\ROp}(\ROp): {}_{\ROp}\E\rightarrow{}_{\ROp}\E$
associated to the operad $\ROp$
considered as an algebra over itself in right modules over itself.
\end{enumerate}

The definition of the functor $\Sym_{\ROp}(N): {}_{\ROp}\E\rightarrow{}_{\POp}\E$
is obviously natural in $N\in{}_{\POp}\M{}_{\ROp}$
so that the map $N\mapsto\Sym_{\ROp}(N)$ defines a functor $\Sym_{\ROp}: {}_{\POp}\M{}_{\ROp}\rightarrow{}_{\POp}\Func{}_{\ROp}$,
where ${}_{\POp}\Func{}_{\ROp}$ denotes the category of functors $F: {}_{\ROp}\E\rightarrow{}_{\POp}\E$
from the category of $\ROp$-algebras in $\E$ to the category of $\POp$-algebras in $\E$.
According to~\cite[Proposition 9.2.1]{FresseModules}:
\begin{enumerate}\setcounter{enumi}{1}
\item\label{Background:FunctorModules:FunctorsToAlgebras:FreeObjects}
For a free $\POp$-algebra in right $\ROp$-modules
we have the identity $\Sym_{\ROp}(\POp(M),A) = \POp(\Sym_{\ROp}(M,A))$,
where on the right hand side we consider the free $\POp$-algebra generated by the object $\Sym_{\ROp}(M,A)\in\E$
associated to $A\in{}_{\ROp}\E$ by the functor $\Sym_{\ROp}(M): {}_{\ROp}\E\rightarrow\E$
determined by $M\in\M_{\ROp}$.
\item\label{Background:FunctorModules:FunctorsToAlgebras:Colimits}
The functor $\Sym_{\ROp}: {}_{\POp}\M{}_{\ROp}\rightarrow{}_{\POp}\Func{}_{\ROp}$ preserves colimits.
\end{enumerate}

\subsection{On extension and restriction of structure}\label{Background:ExtensionRestriction}
Any operad morphism gives rise to adjoint extension and restriction functors
on module categories,
as well as on algebra categories.
The purpose of this subsection is to recall the definition of these functors.

\subsubsection{On extension and restriction of structure for right modules over operads}\label{Background:ExtensionRestriction:RightModules}
On module categories,
the adjoint extension and restriction functors
\begin{equation*}
\psi_!: \M{}_{\ROp}\rightleftarrows\M{}_{\SOp} :\psi^*
\end{equation*}
associated to an operad morphism $\psi: \ROp\rightarrow\SOp$
are very analogous to the classical extension and restriction functors of linear algebra.

The right $\ROp$-module $\psi^* N$ obtained by restriction of structure from an $\SOp$-module $N$
is defined by the object underlying $N$ on which the operad $\ROp$ acts through $\SOp$
by way of the morphism $\psi: \ROp\rightarrow\SOp$.
Usually,
we omit marking the restriction of structure in notation,
unless this abuse of notation creates confusion.

Recall that an operad $\SOp$
forms a bimodule over itself.
By restriction,
we obtain that $\SOp$ is acted on by the operad $\ROp$ on the left
so that $\SOp$ forms an $\ROp$-$\SOp$-bimodule as well.
The extension functor is defined by the relative composition product $\psi_! M = M\circ_{\ROp}\SOp$.

Usually,
we use the expression of the relative composition product $M\circ_{\ROp}\SOp$
to refer to the object $\psi_! M$.
This convention has the advantage of distinguishing extensions of structure on the right from extensions of structure on the left
(whose definition is recalled next, in~\S\ref{Background:ExtensionRestriction:Algebras})
and to stress the analogy with extension of scalars in linear algebra.
Nevertheless we keep using the notation $\psi_!$
to refer to the extension of structure as a functor $\psi_!: \M_{\ROp}\rightarrow\M_{\SOp}$.

The right $\ROp$-module $\psi^* N$
associated to a $\POp$-algebra in right $\SOp$-module $N\in{}_{\POp}\M{}_{\SOp}$,
where $\POp$ is another operad,
inherits an obvious $\POp$-algebra structure and forms a $\POp$-algebra in right $\ROp$-modules.
In the converse direction,
one checks that the relative composition product $M\circ_{\ROp}\SOp$ preserve tensor products,
from which we obtain that the right $\SOp$-module $\psi_! M = M\circ_{\ROp}\SOp$
associated to a $\POp$-algebra in right $\ROp$-modules $M\in{}_{\POp}\M{}_{\ROp}$
inherits a $\POp$-algebra structure
and forms a $\POp$-algebra in right $\SOp$-modules.
Finally,
we have induced extension and restriction functors
\begin{equation*}
\psi_!: {}_{\POp}\M{}_{\ROp}\rightleftarrows{}_{\POp}\M{}_{\SOp} :\psi^*
\end{equation*}
which are obviously adjoint to each other.

\subsubsection{On extension and restriction of structure for algebras over operads}\label{Background:ExtensionRestriction:Algebras}
An operad morphism $\phi: \POp\rightarrow\QOp$ yields adjoint extension and restriction functors
on algebra categories
\begin{equation*}
\phi_!: {}_{\POp}\E\rightleftarrows{}_{\QOp}\E :\phi^*,
\end{equation*}
for any symmetric monoidal category $\E$ over the base category of dg-modules $\C$.

Again,
the $\POp$-algebra $\phi^* B$ obtained by restriction of structure from a $\QOp$-algebra $B$
is defined by the object underlying $B$ on which the operad $\POp$ acts through $\QOp$
by the morphism $\phi: \POp\rightarrow\QOp$.
In the other direction,
the $\POp$-algebra $\phi_! A$ obtained by extension of structure from a $\POp$-algebra $A$
is just characterized by the adjunction relation
\begin{equation*}
\Mor_{{}_{\POp}\E}(\phi_! A,B) = \Mor_{{}_{\QOp}\E}(A,\phi^* B).
\end{equation*}
In fact,
the $\POp$-algebra $\phi_! A$
can be identified with the object $\Sym_{\POp}(\QOp,A)\in\QOp$
associated to $A$ by the functor $\Sym_{\POp}(\QOp): {}_{\POp}\E\rightarrow{}_{\QOp}\E$
where the operad $\QOp$ is considered as an algebra over itself in right modules over $\POp$
(use the restriction of structure on the right of~\S\ref{Background:ExtensionRestriction:RightModules}).

In the case $\E = \M$ and $\E = \M_{\ROp}$,
we obtain extension and restriction functors for left modules over operads
\begin{equation*}
\phi_!: {}_{\POp}\M\rightleftarrows{}_{\QOp}\M :\phi^*
\end{equation*}
and extension and restriction functors on the left for bimodules over operads
\begin{equation*}
\phi_!: {}_{\POp}\M{}_{\ROp}\rightleftarrows{}_{\QOp}\M{}_{\ROp} :\phi^*.
\end{equation*}
In the context of bimodules,
the extension and restriction of structure on the left commute with the extension and restriction of structure on the right.

From the module point of view,
the extension and the restriction of structure on the left
is defined in a symmetric fashion to the extension and the restriction of structure on the right.
In particular,
for extension of structure, we have an identity $\phi_! M = \QOp\circ_{\POp} M$.
Nevertheless we prefer to use the notation $\phi_! M$ to refer to an extension of structure on the left,
rather than the notation of a relative composition product,
because we view the functor $\phi_!: M\mapsto\phi_! M$ as an instance of an extension of structure of $\ROp$-algebras
and this point of view reflects the nature of extensions of structure on the left more properly.

\subsubsection{On extension and restriction of functors}\label{Background:ExtensionRestriction:Functors}
According to~\cite[\S 7.2]{FresseModules},
extension and restriction of structure of modules over operads
reflect extension and restriction operations at the functor level.
For extensions on the right,
we have natural isomorphisms
\begin{equation*}
\Sym_{\SOp}(M\circ_{\ROp}\SOp,B)\simeq\Sym_{\ROp}(M,\psi^* B),
\end{equation*}
for every $M\in\M{}_{\ROp}$ and all $B\in{}_{\SOp}\E$,
as well as natural isomorphisms
\begin{equation*}
\Sym_{\ROp}(N,A)\simeq\Sym_{\SOp}(N,\psi_! A),
\end{equation*}
for every $N\in\M{}_{\SOp}$ and all $A\in{}_{\ROp}\E$,
and similarly in the context of bimodules over operads $M\in{}_{\POp}\M{}_{\ROp}$, $N\in{}_{\POp}\M{}_{\SOp}$
(in this context, the identities hold in the category of $\POp$-algebras).
Symmetrically, for extensions on the left,
we have identities of $\POp$-algebras
\begin{equation*}
\Sym_{\ROp}(\phi_! M,A)\simeq\phi_!\Sym_{\ROp}(M,A),
\end{equation*}
for every $M\in{}_{\POp}\M{}_{\ROp}$,
and
\begin{equation*}
\Sym_{\ROp}(\phi^* N,A)\simeq\phi^*\Sym_{\ROp}(N,A),
\end{equation*}
for every $N\in{}_{\QOp}\M{}_{\ROp}$, where in both cases $A\in{}_{\ROp}\E$.

\setcounter{section}{0}
\renewcommand{\thesection}{\arabic{section}}
\renewcommand{\thesubsection}{\thesection.\arabic{subsection}}
\renewcommand{\thesubsubsection}{\thesubsection.\arabic{subsubsection}}

\mypart{The bar construction and its multiplicative structure}

In this part,
we apply the general theory recalled in \S\S\ref{Background:SymmetricMonoidalCategories}-\ref{Background:ExtensionRestriction}
to prove our main results on the bar construction.

In~\S\ref{BarConstruction},
we recall the definition of the bar construction of differential graded algebras
and we check that this construction is an instance of a functor determined by a module over an operad,
the \emph{bar module}.
For this purpose,
we observe that a generalized bar construction is defined in the setting of modules over operads.
In fact,
the bar module is an instance of a bar construction in that category,
where an operad is considered as an algebra over itself in right modules over itself.
This idea is also used to check homotopical properties of the bar module
associated to an operad.

The multiplicative structure of the bar complex
is examined in~\S\ref{BarStructure},
where we use constructions of~\S\ref{BarConstruction}
to prove the existence and uniqueness of an $E_\infty$-structure on the bar complex
of $E_\infty$-algebras.

In~\S\ref{CategoricalBarConstruction},
we recall the definition of a categorical analogue of the bar construction,
where tensor products are replaced by categorical coproducts,
and we check that this categorical bar construction
forms an instance of a functor determined by a module over an operad.
By~\cite[Theorems 3.2 and 3.5]{Mandell},
the categorical bar construction defines a model of the suspension in the homotopy category of algebras over an operad.
In~\S\ref{HomotopyInterpretation},
we use an equivalence of modules over operads
to prove that the usual bar construction is equivalent to the categorical bar construction
as an $E_\infty$-algebra,
from which we conclude that the usual bar construction
defines a model of the suspension in the homotopy category of $E_\infty$-algebras.
This relationship is used in the next part to deduce from results of~\cite{Mandell}
that the bar complex of a cochain algebra $C^*(X)$
is equivalent as an $E_\infty$-algebra
to $C^*(\Omega X)$,
the cochain algebra of the loop space of $X$.

\subsection*{Conventions}
In the remainder of the article,
the notation $\E$ refers either to the category of dg-modules $\E = \C$
or to a category of right modules over an operad $\E = \M_{\SOp}$
and we do not consider other examples of symmetric monoidal categories
over dg-modules.
The concept of a symmetric monoidal category over dg-modules
is essential to understand our arguments,
but in applications we are only interested in these examples.

From now on,
we use the subcategory $\Op_0\subset\Op$ formed by operads $\POp$ such that $\POp(0) = 0$,
and we assume tacitely that any given operad satisfies this condition.
The assumption $\POp(0) = 0$
amounts to considering algebras without $0$-ary operations $\lambda: \POp(0)\rightarrow A$.
In the sequel,
we say that an operad $\POp$ which has $\POp(0) = 0$ is \emph{non-unitary} and that the associated algebras are \emph{non-unital}.
This setting simplifies the definition of the bar complex
(see~\S\ref{BarConstruction:AinfinityBarComplex:Definition} and~\S\ref{BarConstruction:AinfinityBarComplex:UnitaryContext}).

\section{The bar construction and the bar module}\label{BarConstruction}

\subsection*{Introduction}
In this section,
we check that the bar construction $A\mapsto B(A)$
is identified with the functor associated to a right module over Stasheff's operad
and we check properties of this module.

For our needs,
we study restrictions of the bar construction to categories of algebras over operads $\ROp$,
where $\ROp$ is any operad under Stasheff's operad $\KOp$.
In this context,
we prove:

\begin{mainprop}\label{BarConstruction:prop:BarModuleFeature}
Let $\ROp$ be any operad under Stasheff's operad $\KOp$.
There is a right $\ROp$-module naturally associated to $\ROp$, the bar module $B_{\ROp}$,
such that $B(A) = \Sym_{\ROp}(B_{\ROp},A)$,
for all $A\in{}_{\ROp}\E$.
\end{mainprop}

\medskip
In~\S\ref{BarConstruction:AinfinityBarComplex},
we recall the definition of Stasheff's operad $\KOp$ and the definition of the bar construction
for algebras over this operad.
In~\S\ref{BarConstruction:OperadModuleBarComplex},
we study the bar construction of a $\KOp$-algebra in a category of right modules over an operad $\ROp$.

In~\S\ref{BarConstruction:UnderAinfinityOperads},
we note that all $E_\infty$-operads form operads under $\KOp$.
As a consequence we obtain that all $E_\infty$-algebras have an associated bar complex.
In~\S\ref{BarConstruction:BarModule},
we use the generalized bar complex of $\KOp$-algebras in right modules over an operad $\ROp$
to define the bar module $B_{\ROp}$ associated to an operad $\ROp$ under $\KOp$.
For this aim, we just observe that an operad under $\KOp$ forms a $\KOp$-algebra in right modules over itself.
We study the structure of this right $\ROp$-module $B_{\ROp}$ and the functoriality of the construction $\ROp\mapsto B_{\ROp}$.

\subsection{On Stasheff's operad and the bar complex}\label{BarConstruction:AinfinityBarComplex}
In this section,
we use that the structure of an algebra over an operad $\POp$,
defined by a collection of evaluation morphisms
\begin{equation*}
\lambda: \POp(n)\otimes A^{\otimes n}\rightarrow A,
\end{equation*}
amounts to associating an actual operation $p: A^{\otimes n}\rightarrow A$
to any homogeneous element $p\in\POp(n)$,
at least in the case $\E = \C$, the category of dg-modules,
and $\E = \M_{\ROp}$, the category of right modules over an operad $\ROp$.

For this purpose,
we use the adjunction relation
\begin{equation*}
\Mor_{\E}(\POp(n)\otimes A^{\otimes n},A) = \Mor_{\C}(\POp(n),\Hom_{\E}(A^{\otimes n},A))
\end{equation*}
and an explicit representation of the dg-hom
\begin{equation*}
\Hom_{\E}(-,-): \E^{op}\times\E\rightarrow\C
\end{equation*}
on these categories.

In the context of dg-modules,
an element $f\in\Hom_{\C}(C,D)$ is simply a homogeneous map $f: C\rightarrow D$.
In the context of right modules over an operad,
an element $f\in\Hom_{\M_{\ROp}}(M,N)$
consists of a collection of homogeneous maps of dg-modules $f: M(n)\rightarrow N(n)$, $n\in\NN$,
which commute with the action of symmetric groups
and so that the action of the operad $\ROp$ is preserved by $f: M\rightarrow N$.
In general,
the evaluation morphism of a $\POp$-algebra associates an element $p\in\Hom_{\E}(A^{\otimes n},A)$
to any operation $p\in\POp(n)$.

The standard bar complex is an instance of a construction where the internal differential of a dg-module $C$
is twisted by a cochain $\partial\in\Hom_{\C}(C,C)$
to produce a new dg-module,
which has the same underlying graded module as $C$,
but whose differential is given by the sum $\delta+\partial: C\rightarrow C$.
One has simply to assume that a twisting cochain $\partial$ satisfies the equation $\delta(\partial) + \partial^2 = 0$
in $\Hom_{\C}(C,C)$
to obtain that the map $\delta+\partial$ verifies the equation of differentials $(\delta+\partial)^2 = 0$.
This construction makes sense in the context of right modules over an operad.
In this case,
the twisting cochain $\partial: M\rightarrow M$ is supposed to represents an element of $\Hom_{\M_{\ROp}}(M,M)$
and this condition ensures that the sum $\delta+\partial: M\rightarrow M$
defines a differential of right $\ROp$-modules (for details, compare with definitions of~\cite[\S 2.1.11]{FressePartitions}).

From these observations,
a bar complex in the category of right modules over an operad
can be defined in parallel to the standard bar complex in dg-modules.
Before doing this construction,
we recall the definition of Stasheff's operad,
at least for the sake of completeness.

\subsubsection{On the chain operad of Stasheff's associahedra}\label{BarConstruction:AinfinityBarComplex:StasheffOperad}
The structure of Stasheff's operad $\KOp$ is specified by a pair $\KOp = (\Free(M),\partial)$,
where $\Free(M)$ is a free operad
and $\partial: \Free(M)\rightarrow\Free(M)$ is an operad derivation
that defines the differential of~$\KOp$.
The generating $\Sigma_*$-module $M$ is given by
\begin{equation*}
M(r) = \begin{cases} 0, & \text{if $r = 0,1$}, \\
\Sigma_r\otimes\kk\,\mu_r, & \text{otherwise},
\end{cases}
\end{equation*}
where $\mu_r$ is a generating operation of degree $r-2$.
The derivation $\partial: \Free(M)\rightarrow\Free(M)$ is determined on generating operations
by the formula
\begin{equation*}
\partial(\mu_r) = \sum_{s+t-1 = r}\Bigl\{\sum_{i=1}^{s} \pm\mu_s\circ_i\mu_t\Bigr\}.
\end{equation*}

Let $\AOp$ be operad of associative algebras.
The Stasheff operad is endowed with an operad equivalence $\epsilon: \KOp\xrightarrow{\sim}\AOp$
defined by $\epsilon(\mu_r) = 0$ for $r>2$
and $\epsilon(\mu_2) = \mu$,
where $\mu\in\AOp(2)$ is the operation which represents the product of associative algebras.

\subsubsection{The bar complex}\label{BarConstruction:AinfinityBarComplex:Definition}
Let $A$ be a $\KOp$-algebra in $\E$,
where $\E = \C$, the category of dg-modules,
or $\E = \M_{\ROp}$, the category right modules over an operad.

The (reduced) bar complex of $A$
is defined by the pair $B(A) = (\Tens^c(\Sigma A),\partial)$
formed by the (non-augmented) tensor coalgebra
\begin{equation*}
\Tens^c(\Sigma A) = \bigoplus_{n=1}^{\infty} (\Sigma A)^{\otimes n}
\end{equation*}
where $\Sigma A$ is the suspension of $A$ in $\E$,
together with a twisting cochain
$\partial\in\Hom_{\E}(\Tens^c(\Sigma A),\Tens^c(\Sigma A))$,
called the bar coderivation,
defined pointwise by the formula
\begin{equation*}
\partial(a_1\otimes\dots\otimes a_n)
= \sum_{r=2}^{n}\Bigl\{\sum_{i=1}^{n-r+1} \pm a_1\otimes\dots
\otimes\mu_r(a_i,\dots,a_{i+r-1})\otimes\dots
\otimes a_n\Bigr\}.
\end{equation*}
The internal differential of the bar complex $B(A)$ is the sum $\delta+\partial$ of the natural differential of the tensor coalgebra
$\delta: \Tens^c(\Sigma A)\rightarrow\Tens^c(\Sigma A)$,
induced by the internal differential of $A$,
with the bar coderivation
$\partial: \Tens^c(\Sigma A)\rightarrow\Tens^c(\Sigma A)$,
determined by the $\KOp$-operad action.

In the case of an associative algebra,
the bar coderivation reduces to terms
\begin{equation*}
\partial = \sum_{i=1}^{n-1} \pm a_1\otimes\dots\otimes\mu(a_{i},a_{i+1})\otimes\dots\otimes a_n
\end{equation*}
since the operations $\mu_r\in\KOp(r)$ vanish in $\AOp(r)$ for $r>2$.
Hence, in this case, we recover the standard definition of the bar complex of associative algebras.

According to the definition,
the bar complex forms naturally a coalgebra in $\E$,
but we do not use coalgebra structures further in this article.

\subsubsection{Remark}
In the context of right modules over an operad $\E = \M_{\ROp}$,
we have essentially to form the tensor coalgebra $\Tens^c(\Sigma A)$
in $\M_{\ROp}$.
The suspension of an object $M\in\M_{\ROp}$ can be defined by a tensor product
$\Sigma M = \bar{N}_*(S^1)\otimes M$,
as in the context of dg-modules,
where $\bar{N}_*(S^1)$ is the reduced normalized chain complex of the circle.
The pointwise definition of the bar coderivation $\partial: \Tens^c(\Sigma A)\rightarrow\Tens^c(\Sigma A)$,
makes sense if we recall that $(\Sigma A)^{\otimes n}$
is generated by tensors $a_1\otimes\cdots\otimes a_n\in\Sigma A(r_1)\otimes\cdots\otimes\Sigma A(r_n)$
(we apply the principle of generalized point-tensors of~\cite[\S 0.5]{FresseModules}).
In both cases $\E = \C$ and $\E = \M_{\ROp}$,
the bar coderivation can be defined as a sum of homomorphism tensor products
\begin{equation*}
\id\otimes\dots\otimes\mu_r\otimes\dots\otimes\id\in\Hom_{\E}((\Sigma A)^{\otimes n},(\Sigma A)^{\otimes n-r+1})
\end{equation*}
as well.

\subsubsection{Remark}\label{BarConstruction:AinfinityBarComplex:UnitaryContext}
The definition of~\S\ref{BarConstruction:AinfinityBarComplex:Definition} is the right one for a non-unital algebra.
Similarly,
we consider a non-augmented tensor coalgebra in the definition of $B(A)$,
or equivalently the augmentation ideal of the standard tensor coalgebra,
so that our bar complex forms a non-unital object.
In general it is simpler for us to deal with non-unital algebras and therefore we take this convention.
In the unital context we have to assume that $A$ is augmented
and,
in the definition of $B(A)$,
we have to replace the algebra $A$ by its augmentation ideal $\bar{A}$.

\subsection{The generalized bar complex}\label{BarConstruction:OperadModuleBarComplex}
The bar construction gives by definition a functor $B: {}_{\KOp}\E\rightarrow\E$,
for $\E = \C$ and $\E = \M_{\ROp}$.
In this subsection
we check that standard properties of the usual bar construction of $\KOp$-algebras in dg-modules
hold in the context of right modules over an operad $\ROp$.

First,
we have the easy propositions:

\begin{prop}\label{BarConstruction:OperadModuleComplex:ExtensionRestriction}
Let $\psi: \ROp\rightarrow\SOp$ be any operad morphism.
For any $\KOp$-algebra in right $\ROp$-modules $M$,
we have a natural isomorphism $B(M)\circ_{\ROp}\SOp\simeq B(M\circ_{\ROp}\SOp)$
in the category of right $\SOp$-modules.
\end{prop}

\begin{proof}
Use simply that extension functors $\psi_!(M) = M\circ_{\ROp}\SOp$
commute with tensor products to obtain this isomorphism
(see \cite[\S 7.2]{FresseModules} and recollections in~\S\ref{Background:ExtensionRestriction:RightModules}).
\qed\end{proof}

\begin{prop}\label{BarConstruction:OperadModuleBarComplex:InducedFibrations}
If $\phi: M\rightarrow N$ is a fibration of $\KOp$-algebras in right $\ROp$-modules,
then the induced morphism $B(\phi): B(M)\rightarrow B(N)$
defines a fibration in the category of right $\ROp$-modules
\end{prop}

\begin{proof}
Recall that fibrations in the category of right $\ROp$-modules
are created in the category of dg-modules
and, as such,
are just degreewise epimorphisms.
Therefore
the assertion is an immediate consequence of the definition of the bar complex as a twisted module $B(N) = (\Tens^c(\Sigma N),\partial)$.
Note simply that the tensor coalgebra $\Tens^c(\Sigma N)$ preserves epimorphisms
because the tensor product of right $\ROp$-modules,
inherited from $\Sigma_*$-modules,
has this property.
\qed\end{proof}

Our main task is to check that the bar construction, preserves cofibrations, acyclic cofibrations, and all weak-equivalences
between $\KOp$-algebras which are cofibrant as a right $\ROp$-module.
For this aim we prove that the bar complex has a natural cell decomposition.

Let $D^n$ be the dg-module spanned by an element $e_n$ in degree $n$
and an element $b_{n-1}$ in degree $n-1$
so that $\delta(e_n) = b_{n-1}$.
Consider the submodule $C^{n-1}\subset D^n$
spanned by $b_{n-1}$.
To define the cells,
we use the dg-module embeddings $i_n: C^{n-1}\rightarrow D^n$,
which are generating cofibrations of the category of dg-modules.

\begin{lemm}\label{BarConstruction:OperadModuleBarComplex:CellStructure}
For any $\KOp$-algebra in right $\ROp$-modules $N$,
the bar complex $B(N)$
decomposes into a sequential colimit
\begin{multline*}
0 = B_{\leq 0}(N)\xrightarrow{j_1} B_{\leq 1}(N)\rightarrow\dots
\rightarrow B_{\leq n-1}(N)\xrightarrow{j_n} B_{\leq n}(N)\rightarrow\cdots\\
\dots\rightarrow\colim_n B_{\leq n}(N) = B(N)
\end{multline*}
so that $B_{\leq n}(N)$ is obtained from $B_{\leq n-1}(N)$
by a pushout of the form
\begin{equation*}
\xymatrix{ C^{n-1}\otimes N^{\otimes n}\ar[r]^{f_n}\ar[d]_{i_n} & B_{\leq n-1}(N)\ar@{.>}[d]^{j_n} \\
D^{n}\otimes N^{\otimes n}\ar@{.>}[r]_{g_n} & B_{\leq n}(N) }.
\end{equation*}
This decomposition is also functorial with respect to $N$.
\end{lemm}

\begin{proof}
Indeed,
the object $B(N)$ has a canonical filtration
\begin{equation*}
0 = B_{\leq 0}(N)\hookrightarrow B_{\leq 1}(N)\hookrightarrow\cdots
\hookrightarrow B_{\leq n}(N)\hookrightarrow\cdots
\hookrightarrow\colim_n B_{\leq n}(N) = B(N)
\end{equation*}
defined by
\begin{equation*}
B_{\leq n}(N) = \Tens^c_{\leq n}(\Sigma N) = \bigoplus_{m=1}^{n} (\Sigma N)^{\otimes m}.
\end{equation*}

The summand $(\Sigma N)^{\otimes m}$ is preserved by the natural differential of the tensor coalgebra $\Tens^c(\Sigma N)$.
Moreover the bar coderivation satisfies
\begin{equation*}
\partial((\Sigma N)^{\otimes n})\subset\bigoplus_{r\geq 2} (\Sigma N)^{\otimes n-r+1} = \Tens^c_{\leq n-1}(\Sigma N).
\end{equation*}
Accordingly,
we obtain that $B_{\leq n}(N)$ forms a subobject of $B(N)$
in the category of (differential graded) right $\ROp$-modules.

Besides,
our observation implies that $B_{\leq n}(N)$
splits into a twisted direct sum
\begin{equation*}
B_{\leq n}(N) = (B_{\leq n-1}(N)\oplus(\Sigma N)^{\otimes n},\partial),
\end{equation*}
where $\partial: (\Sigma N)^{\otimes n}\rightarrow B_{\leq n-1}(N)$
represents the restriction of the bar coderivation to the summand $(\Sigma N)^{\otimes n}$
(compare with~\cite[\S 11.2.2]{FresseModules}).
By definition,
the differential of such a twisted object is the sum of the internal differential of $B_{\leq n-1}(N)\oplus(\Sigma N)^{\otimes n}$
with the twisting map $\partial: (\Sigma N)^{\otimes n}\rightarrow B_{\leq n-1}(N)$
on the summand $(\Sigma N)^{\otimes n}$.
Hence the identity $B_{\leq n}(N) = (B_{\leq n-1}(N)\oplus(\Sigma N)^{\otimes n},\partial)$
is obvious.

One checks readily that a twisted direct sum of this form
is equivalent to a pushout of the form of the lemma,
where the attaching map $f_n: C^{n-1}\otimes N^{\otimes n}\rightarrow B_{\leq n-1}(N)$
is yielded by the twisting map $\partial: (\Sigma N)^{\otimes n}\rightarrow B_{\leq n-1}(N)$.
Observe simply that
\begin{equation*}
\Sigma(C^{n-1}\otimes N^{\otimes n}) = \Sigma^{n}(N^{\otimes n}) = (\Sigma N)^{\otimes n}
\end{equation*}
to obtain that any twisting map $\partial: (\Sigma N)^{\otimes n}\rightarrow B_{\leq n-1}(N)$, homogeneous of degree $-1$,
is equivalent to a morphism $f_n: C^{n-1}\otimes N^{\otimes n}\rightarrow B_{\leq n-1}(N)$, homogeneous of degree $0$.
\qed\end{proof}

\begin{prop}\label{BarConstruction:OperadModuleBarComplex:CofibrantStructure}
The bar complex $B(N)$ associated to a $\KOp$-algebra in right $\ROp$-modules $N$ is cofibrant
if the $\KOp$-algebra $N$ defines itself a cofibrant object
in the underlying category of right $\ROp$-modules $\M_{\ROp}$.
\end{prop}

\begin{proof}
The axioms of monoidal model categories
imply that the morphism $i_n\otimes N^{\otimes n}: C^{n-1}\otimes N^{\otimes n}\rightarrow C^{n}\otimes N^{\otimes n}$
forms a cofibration in the category of right $\ROp$-modules $\M_{\ROp}$
if $M$ is cofibrant as a right $\ROp$-module.
As a consequence, we obtain that the morphism $j_n: B_{\leq n-1}(N)\rightarrow B_{\leq n}(N)$
defines a cofibration, for each $n\geq 1$,
since this morphism is obtained by a pushout of $i_n$.
The proposition follows.
\qed\end{proof}

A morphism of $\Sigma_*$-modules $i: M\rightarrow N$
is called a $\C$-cofibration (respectively, an acyclic $\C$-cofibration)
if the morphisms $i: M(n)\rightarrow N(n)$, $n\in\NN$,
are cofibrations in the category of dg-modules $\C$.
Similarly,
a $\Sigma_*$-module $M$ is $\C$-cofibrant
if its underlying collection consists of cofibrant dg-modules.

\begin{lemm}\label{BarConstruction:OperadModuleBarComplex:InducedCofibrations}
Let $i: M\rightarrow N$ be a morphism of $\KOp$-algebras in right $\ROp$-modules
such that the $\KOp$-algebra $M$ is $\C$-cofibrant.

The morphism $B(i): B(M)\rightarrow B(N)$ induced by $i$
forms a $\C$-cofibration (respectively, an acyclic $\C$-cofibration)
if $i$ forms itself a $\C$-cofibration (respectively, an acyclic $\C$-cofibration).
\end{lemm}

\begin{proof}
The morphism $B(i): B(M)\rightarrow B(N)$
can be decomposed naturally into a sequential colimit of morphisms
$j_n: B_{\leq n-1}(N/M)\rightarrow B_{\leq n}(N/M)$,
where
\begin{equation*}
B_{\leq n}(N/M) = B(M)\bigoplus_{B_{\leq n}(M)}B_{\leq n}(N)
\end{equation*}
and $j_n$ is induced componentwise by the morphisms
\begin{equation*}
\xymatrix{ B(M)\ar[d]_{=} & B_{\leq n-1}(M)\ar[d]^{B_{\leq n-1}(i)}\ar[l]\ar[r] & B_{\leq n-1}(N)\ar[d]^{B_{\leq n-1}(i)} \\
B(M) & B_{\leq n}(M)\ar[l]\ar[r] & B_{\leq n}(N) }.
\end{equation*}
One checks readily that $j_n$
fits a pushout of the form
\begin{equation*}
\xymatrix{ {\displaystyle C^{n-1}\otimes N^{\otimes n}\bigoplus_{C^{n-1}\otimes M^{\otimes n}} D^n\otimes M^{\otimes n}}\ar[r]\ar[d] &
{\displaystyle B_{\leq n-1}(N/M)}\ar[d]^{j_n} \\
{\displaystyle D^n\otimes N^{\otimes n}}\ar[r] & {\displaystyle B_{\leq n}(N/M)} }.
\end{equation*}

The underlying dg-modules of the tensor power $M^{\otimes r}$, where $M$ is any right $\ROp$-module,
have an expansion of the form:
\begin{align*}
M^{\otimes r}(m)
& = \bigoplus_{m_1+\cdots+m_r = m} \Sigma_m\otimes_{\Sigma_{m_1}\times\dots\times\Sigma_{m_r}} M(m_1)\otimes\cdots\otimes M(m_r)\\
& = \bigoplus_{m_1+\cdots+m_r = m} (\Sigma_m/\Sigma_{m_1}\times\dots\times\Sigma_{m_r})\otimes M(m_1)\otimes\cdots\otimes M(m_r),
\end{align*}
where the tensor product of the dg-module $T = M(m_1)\otimes\cdots\otimes M(m_r)$
with the coset $K = \Sigma_m/\Sigma_{m_1}\times\dots\times\Sigma_{m_r}$
is defined by a sum of copies of $T$ indexed by $K$, as usual.
By the monoidal model structure of dg-modules,
we obtain that the morphism $i^{\otimes n}: M^{\otimes n}\rightarrow N^{\otimes n}$
forms a $\C$-cofibration (respectively, an acyclic $\C$-cofibration)
if $i$ is so,
as long as $M$ is $\C$-cofibrant.
Under this assumption,
the pushout product-axiom in dg-modules implies that the left-hand side morphism of the pushout above
is a $\C$-cofibration (respectively, an acyclic $\C$-cofibration),
from which we deduce that our morphism $j_n: B_{\leq n-1}(N/M)\rightarrow B_{\leq n}(N/M)$
forms a cofibration (respectively, an acyclic cofibration)
as well.
The conclusion follows.
(Recall that the forgetful functor which maps a right $\ROp$-module $M$
to its underlying collection of dg-modules $\{M(n)\}_{n\in\NN}$
creates all colimits in $\M_{\ROp}$.
Hence we obtain that $\C$-cofibrations and acyclic $\C$-cofibrations
are preserved by pushouts
in the category of right $\ROp$-modules $\M_{\ROp}$.)
\qed\end{proof}

\begin{prop}\label{BarConstruction:OperadModuleBarComplex:HomotopyInvariance}
The morphism $B(i): B(M)\rightarrow B(N)$
induced by a weak-equivalence of $\KOp$-algebras in right $\ROp$-modules $i: M\xrightarrow{\sim} N$
forms itself a weak-equivalence
if the underlying collection of the $\KOp$-algebras $M$ and $N$ consist of cofibrant dg-modules.
\end{prop}

\begin{proof}
Cofibrant algebras over operads form cofibrant objects in the underlying category by~\cite[Corollary 5.5]{BergerMoerdijk}
(see also~\cite[Proposition 12.3.2]{FresseModules}).
This assertion enables us to use the standard Brown's lemma (see for instance~\cite[Lemma 1.1.12]{Hovey})
to obtain the proposition as an immediate consequence of Lemma~\ref{BarConstruction:OperadModuleBarComplex:InducedCofibrations}.
\qed\end{proof}

\subsection{Operads under Stasheff's operad and the bar complex}\label{BarConstruction:UnderAinfinityOperads}
In this subsection,
we examine restrictions of the bar complex to categories of algebras associated to operads $\POp$
equipped with a morphism $\eta: \KOp\rightarrow\POp$.
For our purpose,
we record that any $E_\infty$-operad $\EOp$
can be equipped with such a morphism $\eta: \KOp\rightarrow\EOp$,
so that any algebra over an $E_\infty$-operad
has a bar complex.
By the way,
we recall the definition of an $E_\infty$-operad,
at least to fix conventions.

\subsubsection{Operads under Stasheff's operad and the bar complex}\label{BarConstruction:UnderAinfinityOperads:BarComplexFunctor}
The category of (non-unitary) operads under $\KOp$,
for which we use the notation $\Op_0\backslash\KOp$,
is the comma category of operad morphisms $\eta: \KOp\rightarrow\POp$,
where $\POp\in\Op_0$.
According to this definition,
an operad under Stasheff's operad $\KOp$ is defined by a pair $(\POp,\eta)$
formed by an operad $\POp$ together with an operad morphism $\eta: \KOp\rightarrow\POp$.
Usually,
we omit abusively the morphism $\eta: \KOp\rightarrow\POp$
in the notation of an operad under $\KOp$
and
we identify an object of $\Op_0\backslash\KOp$ with a non-unitary operad $\POp$
endowed with a morphism $\eta: \KOp\rightarrow\POp$
given with $\POp$.

If $\POp$ is an operad under $\KOp$,
then the category of $\POp$-algebras is equipped with a canonical restriction functor $\eta^*: {}_{\POp}\E\rightarrow{}_{\KOp}\E$
associated to the morphism $\eta: \KOp\rightarrow\POp$.
As a consequence,
the bar complex restricts naturally to a functor on the category of $\POp$-algebras,
for all operads $\POp\in\Op_0\backslash\KOp$.
Formally,
this functor is given by the composite
\begin{equation*}
{}_{\POp}\E\xrightarrow{\eta^*}{}_{\KOp}\E\xrightarrow{B}\E.
\end{equation*}
Observations of~\S\ref{BarConstruction:AinfinityBarComplex:Definition}
imply that we recover the usual bar complex of associative algebras
in the case where $\POp$ is the associative operad $\AOp$
together with the canonical augmentation morphism $\epsilon: \KOp\xrightarrow{\sim}\AOp$.

\subsubsection{On $E_\infty$-operads as operads under Stasheff's operad}\label{BarConstruction:UnderAinfinityOperads:EinfinityOperads}
By definition,
an $E_\infty$-operad is an operad $\EOp$
equipped with a weak-equivalence of operads $\epsilon: \EOp\xrightarrow{\sim}\COp$, called the augmentation of $\EOp$,
where $\COp$ denotes the (non-unitary) commutative operad,
the operad associated to the category of (non-unital) associative and commutative algebras.
In the literature,
an $E_\infty$-operad is usually assumed to be $\Sigma_*$-cofibrant
and we take this convention as well.
Observe that the augmentation $\epsilon: \EOp\xrightarrow{\sim}\COp$
is automatically a fibration
because $\COp$ is an operad
in $\kk$-modules, equipped with a trivial differential.

In the introduction of this part, we mention that any $E_\infty$-operad $\EOp$ forms an operad under Stasheff's operad $\KOp$.
Recall that we have an operad morphism $\alpha: \AOp\rightarrow\COp$
so that the restriction functor $\alpha^*: {}_{\COp}\E\rightarrow{}_{\AOp}\E$
represents the embedding from the category of associative and commutative algebras
to the category of all associative algebras.
We simply fix a lifting
\begin{equation*}
\xymatrix{ \KOp\ar@{.>}[r]^{\eta}\ar@{->>}[d]^{\sim} & \EOp\ar@{->>}[d]^{\sim} \\
\AOp\ar[r]^{\alpha} & \COp }
\end{equation*}
in order to obtain an operad morphism $\eta: \KOp\rightarrow\EOp$
such that the restriction functor $\eta^*: {}_{\EOp}\E\rightarrow{}_{\KOp}\E$
extends the standard category embedding $\alpha^*: {}_{\COp}\E\hookrightarrow{}_{\AOp}\E$
from commutative algebras to associative algebras.
Observe that $\eta: \KOp\rightarrow\EOp$ is uniquely determined up to homotopy only.
Therefore,
in this article,
we assume tacitely that such a morphism $\eta: \KOp\rightarrow\EOp$
is fixed for any given $E_\infty$-operad $\EOp$.

By observations of~\S\ref{BarConstruction:UnderAinfinityOperads:BarComplexFunctor},
we obtain that the bar complex restricts to a functor on the category of $\EOp$-algebras.
In addition,
since we have a commutative diagram of restriction functors
\begin{equation*}
\xymatrix{ \E & {}_{\KOp}\E\ar@{.>}[l]_{B} & {}_{\EOp}\E\ar[l]_{\eta^*} \\
& {}_{\AOp}\E\ar[u] & {}_{\COp}\E\ar[l]_{\alpha^*}\ar[u] },
\end{equation*}
we obtain that the bar complex of $\EOp$-algebras extends the usual bar complex
on the category of associative and commutative algebras.

\subsubsection{Remark}
In our construction,
we mention that the morphism $\eta: \KOp\rightarrow\EOp$
is unique up to homotopy.
Indeed,
as usual in a model category,
all morphisms $\eta_0,\eta_1: \KOp\rightarrow\EOp$
that lift the classical operad morphism $\alpha: \AOp\rightarrow\COp$
are connected by a left homotopy in the category of dg-operads.
By~\cite[Theorem 5.2.2]{FresseCylinder},
the existence of such a left homotopy implies the existence of a natural weak-equivalence
between the composite functors
\begin{equation*}
\xymatrix{ {}_{\EOp}\E\ar@<+1mm>[r]^{\eta^*_0}\ar@<-1mm>[r]_{\eta^*_1} & {}_{\KOp}\E\ar[r]^{B} & \E }.
\end{equation*}

To conclude,
we have a well-defined bar complex functor $B: {}_{\EOp}\E\rightarrow\E$
once the $E_\infty$-operad $\EOp$ is provided with a fixed operad morphism $\eta: \KOp\rightarrow\EOp$
that lifts the classical operad morphism $\alpha: \AOp\rightarrow\COp$.
Otherwise
the bar complex functor $B: A\mapsto B(A)$ is uniquely determined up to homotopy only.

\subsection{The bar module}\label{BarConstruction:BarModule}
By definition,
the bar construction of a $\KOp$-algebra in right $\ROp$-modules $N\in{}_{\KOp}\M{}_{\ROp}$
returns a right $\ROp$-module $B(N)$,
and this right $\ROp$-module determines a functor $\Sym_{\ROp}(B(N)): {}_{\ROp}\E\rightarrow\E$.
For our purpose,
we note:

\begin{prop}\label{BarConstruction:AinfinityBarComplex:FunctorBarModule}
Let $\E = \C$, the category of dg-modules, or $\E = \M_{\SOp}$, the category of right modules over an operad $\SOp$.
Let $N$ be any $\KOp$-algebra in right $\ROp$-modules.
The bar complex of $N$ in right $\ROp$-modules satisfies the relation
\begin{equation*}
\Sym_{\ROp}(B(N),A) = B(\Sym_{\ROp}(N,A)),
\end{equation*}
for all $A\in{}_{\ROp}\E$,
where on the right-hand side we consider the bar complex of the $\KOp$-algebra $\Sym_{\ROp}(N,A)\in{}_{\KOp}\E$
associated to $A\in{}_{\ROp}\E$
by the functor $\Sym_{\ROp}(N): {}_{\ROp}\E\rightarrow{}_{\KOp}\E$
defined by $N$.
\end{prop}

\begin{proof}
Since the functor $M\mapsto\Sym_{\ROp}(M)$ preserves internal tensor products of the category of right $\ROp$-modules
and external tensor products over dg-modules,
we obtain
\begin{equation*}
\Sym_{\ROp}(\Tens^c(\Sigma N),A) = \Tens^c(\Sym_{\ROp}(\Sigma N,A)) = \Tens^c(\Sigma \Sym_{\ROp}(N,A)).
\end{equation*}
The map $\partial: \Sym_{\ROp}(\Tens^c(\Sigma N),A)\rightarrow \Sym_{\ROp}(\Tens^c(\Sigma N),A)$
induced by the bar coderivation of $B(N)$
can also be identified with the bar coderivation of $B(\Sym_{\ROp}(N,A))$.
This identification is tautological as the action of $\KOp$ on $\Sym_{\ROp}(N,A)$
is induced by the action of $\KOp$ on $N$ and, hence, the operations $\mu_r: \Sym_{\ROp}(N,A)^{\otimes r}\rightarrow \Sym_{\ROp}(N,A)$
are the maps induced by the operations $\mu_r: N^{\otimes r}\rightarrow N$
on $N$.
\qed\end{proof}

Recall that an operad $\ROp$ forms an algebra over itself in the category of right modules over itself.
If $\ROp$ comes equipped with a morphism $\eta: \KOp\rightarrow\ROp$ and forms an operad under Stasheff's operad $\KOp$,
then $\ROp$ also defines an algebra over $\KOp$ in right modules over itself
by restriction of structure on the left.
The bar module $B_{\ROp}$ is the bar complex $B_{\ROp} = B(N)$ of this $\KOp$-algebra $N = \eta^*\ROp$.
First,
we check that this object fulfils the requirement of Proposition~\ref{BarConstruction:prop:BarModuleFeature}:

\begin{prop}\label{BarConstruction:BarModule:BarModuleFunctor}
Let $\E = \C$, the category of dg-modules, or $\E = \M_{\SOp}$, the category of right modules over an operad $\SOp$.
The functor $\Sym_{\ROp}(B_{\ROp}): {}_{\ROp}\E\rightarrow\E$ associated to the bar module $B_{\ROp}$
is naturally isomorphic to the bar construction $A\mapsto B(A)$
on the category of $\ROp$-algebras in $\E$.
\end{prop}

\begin{proof}
According to Proposition~\ref{BarConstruction:AinfinityBarComplex:FunctorBarModule},
we have $\Sym_{\ROp}(B_{\ROp},A) = B(\Sym_{\ROp}(\eta^*\ROp,A))$,
where $\eta^*\ROp$ is the $\KOp$-algebra in right $\ROp$-modules defined by the operad $\ROp$.
Recall that $\Sym_{\ROp}(\ROp): {}_{\ROp}\E\rightarrow{}_{\ROp}\E$
represents the identity functor of the category of $\ROp$-algebras.
Moreover,
we have an identity $\Sym_{\ROp}(\eta^* N,A) = \eta^*\Sym_{\ROp}(N,A)$
for all $\ROp$-algebras $N$ in right $\ROp$-modules
(see recollections of~\S\ref{Background:ExtensionRestriction:Functors}).
Hence
the object $\Sym_{\ROp}(\eta^*\ROp,A)$ represents the $\KOp$-algebra associated to $A\in{}_{\ROp}\E$
by restriction of structure
and we obtain finally $\Sym_{\ROp}(B_{\ROp},A) = B(\Sym_{\ROp}(\eta^*\ROp,A)) = B(A)$.
\qed\end{proof}

For our purpose,
we examine the functoriality of this construction
with respect to the operad $\ROp$.
For this aim,
we use the following formal observation:

\begin{obsv}\label{BarConstruction:BarModule:Functoriality}
Let $\psi: \ROp\rightarrow\SOp$ be a morphism of operads under $\KOp$.

\begin{enumerate}
\item
The map $\psi: \ROp\rightarrow\SOp$ defines a morphism $\psi_{\sharp}: \ROp\rightarrow\SOp$
in the category of $\KOp$-algebras in right $\ROp$-modules,
where we use
restrictions of structure on the left to make $\ROp$ (respectively, $\SOp$) into a $\KOp$-algebra
and
restrictions of structure on the right to make $\SOp$ into a right $\ROp$-module.
\item
The morphism of $\KOp$-algebras in right $\SOp$-modules $\psi_{\flat}: \ROp\circ_{\ROp}\SOp\rightarrow\SOp$
adjoint to $\psi_{\sharp}: \ROp\rightarrow\SOp$
forms an isomorphism.
\end{enumerate}
\end{obsv}

From this observation and Observation~\ref{BarConstruction:OperadModuleComplex:ExtensionRestriction},
we deduce that a morphism of operads under $\KOp$
gives rise to a morphism
$\psi_{\sharp}: B_{\ROp}\rightarrow B_{\SOp}$,
in the category of right $\ROp$-modules
and to an isomorphism
$\psi_{\flat}: B_{\ROp}\circ_{\ROp}\SOp\xrightarrow{\simeq} B_{\SOp}$,
which is obviously adjoint to $\psi_{\sharp}$.
Since we assume that weak-equivalences (respectively, fibrations) are created by forgetful functors,
we obtain that $\psi_{\sharp}: \ROp\rightarrow\SOp$ defines a weak-equivalence (respectively, a fibration)
in the category of $\KOp$-algebras in right $\ROp$-modules
if $\psi$ is a weak-equivalence (respectively, a fibration) of operads.
Hence,
Proposition~\ref{BarConstruction:OperadModuleBarComplex:InducedFibrations}
and
Proposition~\ref{BarConstruction:OperadModuleBarComplex:HomotopyInvariance}
return:

\begin{prop}\label{BarConstruction:BarModule:BarModuleHomotopy}
The morphism $\psi_{\sharp}: B_{\ROp}\rightarrow B_{\SOp}$
defines a fibration in the category of right $\ROp$-modules
if $\psi: \ROp\rightarrow\SOp$ is a fibration of operads under $\KOp$.

The morphism $\psi_{\sharp}: B_{\ROp}\rightarrow B_{\SOp}$
defines a weak-equivalence in the category of right $\ROp$-modules
if $\psi: \ROp\rightarrow\SOp$ is a weak-equivalence of operads under $\KOp$
and the underlying collections of the operads $\ROp$ and $\SOp$
consist of cofibrant dg-modules $\ROp(n),\SOp(n)\in\C$, $n\in\NN$.
\qed
\end{prop}

The isomorphism
$\psi_{\flat}: B_{\ROp}\circ_{\ROp}\SOp\xrightarrow{\simeq}B_{\SOp}$
has a natural interpretation at the functor level.
In~\S\ref{Background:ExtensionRestriction:Functors},
we recall that the functor $\Sym_{\SOp}(M\circ_{\ROp}\SOp): {}_{\SOp}\E\rightarrow\E$,
where $M\circ_{\ROp}\SOp$ is the extension of structure of a right $\ROp$-module $M$,
is isomorphic to the composite
\begin{equation*}
{}_{\SOp}\E\xrightarrow{\psi^*}{}_{\ROp}\E\xrightarrow{\Sym_{\ROp}(M)}\E,
\end{equation*}
where $\psi^*: {}_{\SOp}\E\rightarrow{}_{\ROp}\E$
is the restriction functor associated to $\psi: \ROp\rightarrow\SOp$.
For an operad under Stasheff's operad,
the bar complex functor $B: {}_{\ROp}\E\rightarrow\E$
is defined precisely by a composite
of this form:
\begin{equation*}
{}_{\ROp}\E\xrightarrow{\eta^*}{}_{\KOp}\E\xrightarrow{B}\E,
\end{equation*}
where we assume again $\E = \C$ or $\E = \M_{\SOp}$.
Now suppose given a diagram
\begin{equation*}
\xymatrix{ & \KOp\ar[dl]_{\eta}\ar[dr]^{\theta} & \\
\ROp\ar[rr]_{\psi} && \SOp }
\end{equation*}
so that $\psi: \ROp\rightarrow\SOp$ is a morphism of operads under $\KOp$.
The diagram of functors
\begin{equation*}
\xymatrix{ & \E & \\ {}_{\ROp}\E\ar[ur]^{B = \Sym_{\ROp}(B_{\ROp})} && {}_{\SOp}\E\ar[ul]_{B = \Sym_{\SOp}(B_{\SOp})}\ar[ll]^{\psi^*} }
\end{equation*}
commutes
just because the relation $\theta = \psi\eta$ implies that the diagram of restriction functors
\begin{equation*}
\xymatrix{ & {}_{\KOp}\E & \\
{}_{\ROp}\E\ar[ur]^{\eta^*} && {}_{\SOp}\E\ar[ul]_{\theta^*}\ar[ll]^{\psi^*} }
\end{equation*}
commutes.
Thus,
for a morphism $\psi: \ROp\rightarrow\SOp$ in $\Op_0\backslash\KOp$,
we have a natural isomorphism
$\Sym_{\ROp}(B_{\ROp},\psi^* A)\simeq \Sym_{\SOp}(B_{\SOp},A)$, for all $A\in{}_{\ROp}\E$.
Moreover:

\begin{prop}\label{BarConstruction:BarModule:Extension}
Let $\psi: \ROp\rightarrow\SOp$ be any morphism of operads under $\KOp$.
The natural isomorphism
\begin{equation*}
\psi_{\flat}: B_{\ROp}\circ_{\ROp}\SOp\xrightarrow{\simeq} B_{\SOp}
\end{equation*}
induces an isomorphism of functors $\Sym_{\SOp}(\psi_*): \Sym_{\SOp}(B_{\ROp}\circ_{\ROp}\SOp)\rightarrow\Sym_{\SOp}(B_{\SOp})$
that fits a commutative diagram
\begin{equation*}
\xymatrix{ \Sym_{\SOp}(B_{\ROp}\circ_{\ROp}\SOp,A)\ar@{.>}[dr]_{\Sym_{\SOp}(\psi_{\flat},A)}\ar[rr]^{\simeq} && \Sym_{\ROp}(B_{\ROp},\psi^* A)\ar[dl]^{\simeq} \\
& \Sym_{\SOp}(B_{\SOp},A) & },
\end{equation*}
for all $A\in{}_{\SOp}\E$.
\end{prop}

\begin{proof}
The proposition
is a formal consequence of coherence properties between distribution isomorphisms
$(M\otimes N)\circ_{\ROp}\SOp\simeq(M\circ_{\ROp}\SOp)\otimes(N\circ_{\ROp}\SOp)$
and the functor isomorphisms $\Sym_{\ROp}(M\otimes N)\simeq \Sym_{\ROp}(M)\otimes \Sym_{\ROp}(N)$.
\qed\end{proof}

In particular,
for the initial morphism $\eta: \KOp\rightarrow\ROp$ of an operad $\ROp\in\Op_0\backslash\KOp$,
the isomorphism $\eta_{\flat}: B_{\KOp}\circ_{\KOp}\ROp\simeq B_{\ROp}$
reflects the definition of $\Sym_{\ROp}(B_{\ROp}): {}_{\ROp}\E\rightarrow\E$
as the restriction of a functor $B: {}_{\KOp}\E\rightarrow\E$.

\medskip
To complete our results,
observe that the operad $\ROp$ defines a cofibrant object in the category of right modules over itself.
Accordingly,
Proposition~\ref{BarConstruction:OperadModuleBarComplex:CofibrantStructure}
implies:

\begin{prop}\label{BarConstruction:BarModule:CofibrantStructure}
The module $B_{\ROp}$ forms a cofibrant object in the category of right $\ROp$-modules.\qed
\end{prop}

\section{The multiplicative structure of the bar construction}\label{BarStructure}

\subsection*{Introduction}
In this section,
we prove the existence and uniqueness of algebra structures on the bar module of $E_\infty$-operads.
Then we use the correspondence between right modules and functors
to obtain the existence and uniqueness of functorial algebra structures
on the bar construction itself $B(A)$,
for all algebras $A$ over a given $E_\infty$-operad $\EOp$.

\medskip
To prove the existence of algebra structures on the bar module $B_{\EOp}$
the idea is to use endomorphism operads of right modules over operads.
Recall briefly that the endomorphism operad of an object $M$ in a category $\E$
is a universal operad in dg-modules $\End_{M}$
such that the structure of a $\POp$-algebra on $M$
is equivalent to an operad morphism $\nabla: \POp\rightarrow\End_M$.
In this section,
we may specify $\POp$-algebra structures by pairs $(M,\nabla)$, where $\nabla: \POp\rightarrow\End_M$ is the operad morphism
that determines the $\POp$-algebra structure of $M$,
because we deal with objects which are not endowed with a natural internal $\POp$-algebra
structure.

In~\S\ref{BarStructure:Existence:CommutativeStructure},
we observe that the bar module $B_{\COp}$ of the commutative operad $\COp$
can be equipped with the structure of a commutative algebra,
like the bar complex of any commutative algebra.
This structure is represented by an operad morphism $\nabla_c: \COp\rightarrow\End_{B_{\COp}}$,
from the commutative operad $\COp$
to the endomorphism operad of~$B_{\COp}$.

Our main existence theorem, proved in~\S\ref{BarStructure:Existence},
reads:

\begin{mainthm}\label{BarStructure:thm:Existence}
Let $\EOp$ be any $E_\infty$-operad.
Let $\QOp$ be any cofibrant operad augmented over the commutative operad $\COp$.
Let $\epsilon: \EOp\rightarrow\COp$ and $\phi: \QOp\rightarrow\COp$
denote the respective augmentations of these operads.

There is an operad morphism $\nabla_{\epsilon}: \QOp\rightarrow\End_{B_{\EOp}}$
which equips the bar module $B_{\EOp}$
with a left $\QOp$-action
so that:
\begin{enumerate}
\item\label{thm:Existence:ModuleAlgebra}
The bar module $B_{\EOp}$ forms a $\QOp$-algebra in right $\EOp$-modules.
\item\label{thm:Existence:CommutativeReduction}
The natural isomorphism of right $\COp$-modules $B_{\EOp}\circ_{\EOp}\COp\simeq B_{\COp}$
defines an isomorphism in the category of $\QOp$-algebras in right $\COp$-modules
\begin{equation*}
(B_{\EOp},\nabla_{\epsilon})\circ_{\EOp}\COp\simeq\phi^*(B_{\COp},\nabla_c),
\end{equation*}
where the $\QOp$-algebra structure of $B_{\COp}$ is obtained by restriction of its $\COp$-algebra structure
through the augmentation of $\QOp$.
\end{enumerate}
\end{mainthm}

The interpretation of this theorem at the level of the bar construction
is straightforward and is also established in~\S\ref{BarStructure:Existence}.

To prove Theorem~\ref{BarStructure:thm:Existence},
we observe that condition~(\ref{thm:Existence:CommutativeReduction})
is equivalent to a lifting problem in the category of operads,
for which axioms of model categories
imply immediately the existence of a solution.

\medskip
In~\S\ref{BarStructure:Uniqueness},
we check that the isomorphism $(B_{\EOp},\nabla_{\epsilon})\circ_{\EOp}\COp\simeq\phi^*(B_{\COp},\nabla_c)$
of condition~(\ref{thm:Existence:CommutativeReduction})
implies, by adjunction,
the existence of a weak-equivalence in the category of $\QOp$-algebras in right $\EOp$-modules
\begin{equation*}
(B_{\EOp},\nabla_{\epsilon})\xrightarrow{\sim}\phi^*(B_{\COp},\nabla_c).
\end{equation*}
From this assertion
we conclude immediately that all solutions of the existence Theorem~\ref{BarStructure:thm:Existence}
yield equivalent objects
in the homotopy category of $\QOp$-algebras in right $\EOp$-modules.
Then
we apply the homotopy invariance theorems of~\cite[\S 15]{FresseModules}
to obtain that all solutions of the existence Theorem~\ref{BarStructure:thm:Existence}
give homotopy equivalent structures on the bar construction.
This gives our uniqueness result.

\numberwithin{mainthm}{subsection}
\renewcommand{\themainthm}{\thesubsection.\Alph{mainthm}}

\subsection{The existence theorem}\label{BarStructure:Existence}
This subsection is devoted to the existence part of our theorems.

To begin with,
we examine the structure of the bar construction of commutative algebras
in the context of right modules over operads
-- we check that the bar module $B_{\COp}$ of the commutative operad $\COp$
forms naturally a commutative algebra in right $\COp$-modules.
Then
we describe constructions which make Theorem~\ref{BarStructure:thm:Existence}
equivalent to a lifting problem in the category of operads
and we solve this lifting problem by arguments of homotopical algebra.

\begin{recollection}[The shuffle product]\label{BarStructure:Existence:ShuffleProduct}
The tensor coalgebra $\Tens^c(\Sigma A)$
can be equipped with a product
$\smile: \Tens^c(\Sigma A)\otimes \Tens^c(\Sigma A)\rightarrow \Tens^c(\Sigma A)$
defined componentwise by sums of tensor permutations
\begin{equation*}
(\Sigma A)^{\otimes m}\otimes(\Sigma A)^{\otimes n}\xrightarrow{\sum_w w_*}(\Sigma A)^{\otimes m+n}
\end{equation*}
where $w$ ranges over the set of $(m,n)$-shuffles in $\Sigma_{m+n}$.
This product is naturally associative and commutative.
For an associative and commutative algebra $A$,
the bar coderivation $\partial: \Tens^c(\Sigma A)\rightarrow \Tens^c(\Sigma A)$
defines a derivation with respect to $\smile$.
Hence, in this case, we obtain that the bar complex $B(A) = (\Tens^c(\Sigma A),\partial)$
is still an associative and commutative algebra.

Clearly,
this standard construction for commutative algebras in dg-modules
can be extended to algebras in a category of right modules
over an operad $\ROp$ --
just use the symmetry isomorphism of the tensor product of right $\ROp$-modules
in the definition of the shuffle product.
Then
we obtain that the bar complex $B(N)$
comes equipped with the structure of a commutative algebra in right $\ROp$-modules
if $N$ is so.
\end{recollection}

Recall that the map $\Sym_{\ROp}(N): A\mapsto \Sym_{\ROp}(N,A)$ defines a functor from $\ROp$-algebras to commutative algebras
if $N$ is a commutative algebra in right $\ROp$-modules.
As the map $\Sym_{\ROp}: M\mapsto \Sym_{\ROp}(M)$ defines a functor of symmetric monoidal categories
(check recollections in~\S\ref{Background:FunctorModules:Operations}),
we obtain that the shuffle product
\begin{equation*}
(\Sigma N)^{\otimes m}\otimes(\Sigma N)^{\otimes n}\xrightarrow{\sum_w w_*}(\Sigma N)^{\otimes m+n}
\end{equation*}
corresponds to the shuffle product
\begin{equation*}
\Sym_{\ROp}(\Sigma N,A)^{\otimes m}\otimes \Sym_{\ROp}(\Sigma N,A)^{\otimes n}\xrightarrow{\sum_w w_*} \Sym_{\ROp}(\Sigma N,A)^{\otimes m+n}
\end{equation*}
at the functor level.
Hence
we obtain finally:

\begin{obsv}\label{BarStructure:Existence:CommutativeStructure}
Let $\E$ be the category of dg-modules $\E = \C$,
or any category of right modules over an operad $\E = \M_{\SOp}$.
For a commutative algebra in right $\ROp$-modules $N$,
the bar complex $B(N)$ comes equipped with the structure of a commutative algebra in right $\ROp$-modules
so that the isomorphism of functors $\Sym_{\ROp}(B(N),A)\simeq B(\Sym_{\ROp}(N,A))$
defines an isomorphism of commutative algebras,
for all $A\in{}_{\ROp}\E$.
\end{obsv}

We apply this observation to the commutative algebra in right $\COp$-modules
formed by the commutative operad itself,
for which we have $\Sym_{\COp}(\COp) = \Id$,
the identity functor on the category of commutative algebras.
We obtain that the standard commutative algebra structure of the bar construction
is realized by the structure of a commutative algebra in right $\COp$-modules
on the bar module $B_{\COp}$.

Recall that,
for any morphism $\psi: \ROp\rightarrow\SOp$ in $\Op_0\backslash\KOp$,
we have a natural isomorphism $\psi_{\flat}: B_{\ROp}\circ_{\ROp}\SOp\simeq B_{\SOp}$
and this relation reflects the definition of the bar complex $B: {}_{\ROp}\E\rightarrow\E$
by the restriction of a functor $B: {}_{\KOp}\E\rightarrow\E$ (see Proposition~\ref{BarConstruction:BarModule:Extension}).
In particular,
for an $E_\infty$-operad $\EOp$, equipped with an augmentation morphism $\epsilon: \EOp\xrightarrow{\sim}\COp$,
we have an isomorphism $\epsilon_{\flat}: B_{\EOp}\circ_{\EOp}\COp\simeq B_{\COp}$.

Our aim is to lift the structure of the bar module $B_{\COp}$ of the commutative operad~$\COp$
to the bar module $B_{\EOp}$ of any $E_\infty$-operad $\EOp$.
Before proving our result,
we recall the definition of an endomorphism operad
and we give an interpretation of extensions and restrictions of structure
in terms of morphisms on endomorphism operads.

\begin{recollection}[Endomorphism operads]\label{BarStructure:Existence:EndomorphismOperad}
The endomorphism operad of an object $M$ in a category $\E$
is defined by the hom-objects
\begin{equation*}
\End_M(n) = \Hom_{\E}(M^{\otimes n},M),
\end{equation*}
where the symmetric groups operate by tensor permutations on the source
and the operad structure of $\End_M$
is deduced from the composition operation of enriched symmetric monoidal categories.
For a $\POp$-algebra $A$,
the operad morphism $\nabla: \POp\rightarrow\End_A$,
equivalent to the $\POp$-algebra structure of $A$,
is defined simply by the morphisms
\begin{equation*}
\nabla: \POp(n)\rightarrow\Hom_{\E}(A^{\otimes n},A)
\end{equation*}
adjoint to the evaluation morphisms $\lambda: \POp(n)\otimes A^{\otimes n}\rightarrow A$.
We refer to~\cite{KapranovManin} or to~\cite[\S 3.4,\S 6.3]{FresseModules}
for an explicit definition of $\End_M$ in the context of right modules over operads.
In the sequel,
we only use general properties of $\End_M$
arising from the abstract definition of hom-objects $\Hom_{\E}(-,-)$.

Note that endomorphism operads $\End_M$
have a $0$-term
\begin{equation*}
\End_{M}(0) = M\not=0
\end{equation*}
in contrast to our conventions on operads, but this apparent contradiction does not create any difficulty:
in our constructions,
one can replace any endomorphism operad $\End_{M}$
by a suboperad $\overline{\End}_{M}\in\Op_0$
such that
\begin{equation*}
\overline{\End}_{M}(n) = \begin{cases} 0, & \text{if $n = 0$}, \\ \End_{M}(n), & \text{otherwise}, \end{cases}
\end{equation*}
because any operad morphism $\nabla: \POp\rightarrow\End_{M}$, where $\POp\in\Op_0$,
factors through~$\overline{\End}_{M}$.
\end{recollection}

\begin{recollection}[Endomorphism operads and extension functors]\label{BarStructure:Existence:EndomorphismExtensionRestriction}
Recall that an operad morphism $\psi: \ROp\rightarrow\SOp$,
gives rise to a functor of extension of structure on the right
\begin{equation*}
\psi_!: {}_{\POp}\M{}_{\ROp}\rightarrow{}_{\POp}\M{}_{\SOp}
\end{equation*}
(see recollections in~\S\ref{Background:ExtensionRestriction}).
One can also observe that the operad morphism $\psi: \ROp\rightarrow\SOp$
induces a morphism
of endomorphism operads:
\begin{equation*}
\psi_!: \End_M\rightarrow\End_{M\circ_{\ROp}\SOp},
\end{equation*}
for all $M\in\M_{\ROp}$,
essentially because the extension functor $\psi_!: M\mapsto M\circ_{\ROp}\SOp$
preserves tensor products
(see~\cite[Proposition 9.4.4]{FresseModules} and recollections in~\S\ref{Background:ExtensionRestriction}).
For a $\POp$-algebra in right $\ROp$-modules represented by a pair $(N,\nabla)$,
where $N\in\M{}_{\ROp}$ and $\nabla: \POp\rightarrow\End_N$,
we obtain that the $\POp$-algebra in right $\SOp$-modules $\psi_!(N,\nabla) = (N,\nabla)\circ_{\ROp}\SOp$,
obtained from $(N,\nabla)$ by extension of structure on the right,
is represented by the pair $\psi_!(N,\nabla) = (\psi_! N,\psi_!\nabla) = (N\circ_{\ROp}\SOp,\nabla\circ_{\ROp}\SOp)$,
where $\psi_!\nabla = \nabla\circ_{\ROp}\SOp$
is the composite
\begin{equation*}
\POp\xrightarrow{\nabla}\End_N\xrightarrow{\psi_!}\End_{N\circ_{\ROp}\SOp}.
\end{equation*}
This assertion is proved by a formal verification (we refer to~\cite[\S 3.4,\S 9.4]{FresseModules}).

The functor of restriction
of structure on the left
\begin{equation*}
\phi^*: {}_{\QOp}\M{}_{\ROp}\rightarrow{}_{\POp}\M{}_{\ROp},
\end{equation*}
where $\phi: \POp\rightarrow\QOp$ is an operad morphism,
has an obvious simpler description
in terms of operad morphisms.
Namely,
for any $\QOp$-algebra in right $\ROp$-modules represented by a pair $(N,\nabla)$,
where $N\in\M{}_{\ROp}$ and $\nabla: \QOp\rightarrow\End_N$,
the $\POp$-algebra in right $\ROp$-modules $\phi^*(N,\nabla)$,
obtained from $(N,\nabla)$ by restriction of structure on the left,
is represented by the pair $\phi^*(N,\nabla) = (N,\nabla\phi)$,
where $\nabla\phi$ is the composite
\begin{equation*}
\POp\xrightarrow{\phi}\QOp\xrightarrow{\nabla}\End_{N}.
\end{equation*}
\end{recollection}

For bar modules,
we have an isomorphism $\psi_{\flat}: B_{\ROp}\circ_{\ROp}\SOp\xrightarrow{\simeq} B_{\SOp}$
and hence an isomorphism of endomorphism operads
\begin{equation*}
\End_{B_{\ROp}\circ_{\ROp}\SOp}\simeq\End_{B_{\SOp}}.
\end{equation*}
Accordingly,
we obtain that any morphism $\psi: \ROp\rightarrow\SOp$
in the category $\Op_0\backslash\KOp$ of operads under Stasheff's operad $\KOp$
gives rise to a morphism
\begin{equation*}
\End_{B_{\ROp}}\xrightarrow{\psi_*}\End_{B_{\SOp}}.
\end{equation*}
One checks readily that the map $\psi\mapsto\psi_*$ preserves composites and identities
so that the map $\ROp\mapsto\End_{B_{\ROp}}$
defines a functor on $\Op_0\backslash\KOp$.
In addition,
we have:

\begin{obsv}\label{BarStructure:Existence:LiftingProblem}
Let $\EOp$ be any $E_\infty$-operad, equipped with an augmentation $\epsilon: \EOp\xrightarrow{\sim}\COp$.
Let $\QOp$ be any operad together with an augmentation $\phi: \QOp\rightarrow\COp$.

Let $\nabla_c: \COp\rightarrow\End_{B_{\COp}}$
be the morphism of dg-operads
determined by the commutative algebra structure of the bar module $B_{\COp}$.
Let $\nabla_{\epsilon}: \QOp\rightarrow\End_{B_{\EOp}}$
be an operad morphism
which provides the bar module $B_{\EOp}$ with the structure of a $\QOp$-algebra in right $\EOp$-modules.

The natural isomorphism $\epsilon_{\flat}: B_{\EOp}\circ_{\EOp}\COp\xrightarrow{\simeq} B_{\COp}$
defines an isomorphism
in the category of $\QOp$-algebras in right $\COp$-modules
\begin{equation*}
\epsilon_{\flat}: (B_{\EOp},\nabla_{\epsilon})\circ_{\EOp}\COp\xrightarrow{\simeq}\phi^*(B_{\COp},\nabla_c)
\end{equation*}
if and only if $\nabla_{\epsilon}$ fits a commutative diagram
\begin{equation*}
\xymatrix{ \QOp\ar[d]_{\phi}\ar@{.>}[r]^(0.4){\nabla_{\epsilon}} & \End_{B_{\EOp}}\ar[d]^{\epsilon_*} \\ \COp\ar[r]_(0.4){\nabla_c} & \End_{B_{\COp}} }.
\end{equation*}
\end{obsv}

To solve the lifting problem arising from this assertion,
we prove:

\begin{lemm}\label{BarStructure:Existence:EndomorphismOperadFibrations}
The functor $\ROp\mapsto\End_{B_{\ROp}}$ preserves fibrations and acyclic fibrations
between operads $\ROp\in\Op_0\backslash\KOp$ whose underlying collection $\ROp(n)$, $n\in\NN$,
consists of cofibrant dg-modules.
\end{lemm}

\begin{proof}
In this proof,
we prefer to use the notation of the functor $\psi_!: \M_{\ROp}\rightarrow\M_{\SOp}$
to denote extensions of structure of right modules over operads
rather than the equivalent relative composition product $\psi_! M = M\circ_{\ROp}\SOp$.
Similarly,
we use the notation of the functor $\psi^*: \M_{\ROp}\rightarrow\M_{\SOp}$
to denote the restriction of structure of right modules over operads.

In general,
the morphism of endomorphism operads $\psi_!: \End_{M}\rightarrow\End_{\psi_! M}$ induced by an operad morphism $\psi: \ROp\rightarrow\SOp$
consists of morphisms
\begin{equation*}
\Hom_{\ROp}(M^{\otimes r},M)\xrightarrow{\psi_!}\Hom_{\SOp}(\psi_! M^{\otimes r},\psi_! M)
\end{equation*}
formed by using that $\psi_!: \M{}_{\ROp}\rightarrow\M{}_{\SOp}$
defines a functor of symmetric monoidal categories over dg-modules.
We use the adjunction between extension and restriction functors $\psi_!: \M{}_{\ROp}\rightleftarrows\M{}_{\SOp} :\psi^*$
to identify these morphisms with composites
\begin{equation*}
\Hom_{\ROp}(M^{\otimes r},M)
\xrightarrow{\eta(M)_*}\Hom_{\ROp}(M^{\otimes r},\psi^*\psi_! M)
\xrightarrow{\simeq}\Hom_{\SOp}(\psi_! M^{\otimes r},\psi_! M),
\end{equation*}
where $\eta(M)_*$ refers to the morphism on hom-objects induced by the adjunction unit
$\eta(M): M\rightarrow\psi^*\psi_!(M)$.

For endomorphism operads of bar modules,
we obtain that the morphism
$\psi_*: \End_{B_{\ROp}}\rightarrow\End_{B_{\SOp}}$,
induced by a morphism $\psi: \ROp\rightarrow\SOp$ in $\Op_0\backslash\KOp$, can be defined by composites
in diagrams of the form:
\begin{equation*}
\xymatrix{ \Hom_{\ROp}(B_{\ROp}^{\otimes r},B_{\ROp})\ar[r]^(0.45){\eta(B_{\ROp})_*}\ar@/_/[dr]!L_{(\psi_{\sharp})_*} &
\Hom_{\ROp}(B_{\ROp}^{\otimes r},\psi^*\psi_! B_{\ROp})\ar[r]^{\simeq}\ar[d]^{\simeq}_{\psi^*(\psi_{\flat})_*} &
\Hom_{\SOp}(\psi_! B_{\ROp}^{\otimes r},\psi_! B_{\ROp})\ar[d]^{\simeq} \\
&
\Hom_{\ROp}(B_{\ROp}^{\otimes r},\psi^* B_{\SOp})\ar[r]_{\simeq} &
\Hom_{\SOp}(\psi_! B_{\ROp}^{\otimes r},B_{\SOp}) \\
&& \Hom_{\SOp}(B_{\SOp}^{\otimes r},B_{\SOp})\ar[u]_{\simeq} }.
\end{equation*}
Accordingly,
to prove our lemma,
we are reduced to check that the morphism
\begin{equation*}
\Hom_{\ROp}(B_{\ROp}^{\otimes r},B_{\ROp})\xrightarrow{(\psi_{\sharp})_*}
\Hom_{\ROp}(B_{\ROp}^{\otimes r},\psi^* B_{\SOp})
\end{equation*}
induced by $\psi_{\sharp}: B_{\ROp}\rightarrow\psi^* B_{\SOp}$
forms a fibration (respectively, an acyclic fibration) of dg-modules
if $\psi: \ROp\rightarrow\SOp$
is so.

In Lemma~\ref{BarConstruction:BarModule:CofibrantStructure},
we prove that $B_{\ROp}$ forms a cofibrant right $\ROp$-module.
In Proposition~\ref{BarConstruction:BarModule:BarModuleHomotopy},
we record that $\psi_{\sharp}: B_{\ROp}\rightarrow\psi^* B_{\SOp}$
forms a fibration (respectively, an acyclic fibration) of right $\ROp$-modules
if $\psi$ is a fibration (respectively, an acyclic fibration) of operads.
By axioms of symmetric monoidal model categories enriched over dg-modules,
we can conclude from these assertions that the morphism
\begin{equation*}
\Hom_{\ROp}(B_{\ROp}^{\otimes r},B_{\ROp})\xrightarrow{(\psi_{\sharp})_*}
\Hom_{\ROp}(B_{\ROp}^{\otimes r},\psi^* B_{\SOp})
\end{equation*}
forms a fibration (respectively, an acyclic fibration)
if $\psi: \ROp\rightarrow\SOp$
is so
and this proves the lemma.
\qed\end{proof}

By axioms of model categories,
Lemma~\ref{BarStructure:Existence:EndomorphismOperadFibrations}
implies immediately:

\begin{lemm}\label{BarStructure:Existence:LiftingProblemSolutions}
Let $\EOp$ be any $E_\infty$-operad, equipped with an augmentation $\epsilon: \EOp\xrightarrow{\sim}\COp$.
Let $\QOp$ be any operad together with an augmentation $\phi: \QOp\rightarrow\COp$.
If $\QOp$ is cofibrant,
then the lifting problem
\begin{equation*}
\xymatrix{ \QOp\ar[d]_{\phi}\ar@{.>}[r]^(0.4){\nabla_{\epsilon}} & \End_{B_{\EOp}}\ar[d]^{\epsilon_*} \\ \COp\ar[r]_(0.4){\nabla_c} & \End_{B_{\COp}} }
\end{equation*}
has a solution.\qed
\end{lemm}

From which we conclude:

\begin{mainthm}[{Claim of Theorem~\ref{BarStructure:thm:Existence}}]\label{BarStructure:Existence:ModuleLevel}
Assume that $\QOp$ is cofibrant.
Then there is a morphism $\nabla_{\epsilon}: \QOp\rightarrow\End_{B_{\EOp}}$,
which provides the bar module $B_{\EOp}$ with the structure of a $\QOp$-algebra in right $\EOp$-modules,
and so that the natural isomorphism of right $\COp$-modules $\epsilon_{\flat}: B_{\EOp}\circ_{\EOp}\COp\xrightarrow{\simeq} B_{\COp}$
defines an isomorphism
\begin{equation*}
\epsilon_{\flat}: (B_{\EOp},\nabla_{\epsilon})\circ_{\EOp}\COp\xrightarrow{\simeq}\phi^*(B_{\COp},\nabla_c)
\end{equation*}
in the category of $\QOp$-algebras in right $\COp$-modules.\qed
\end{mainthm}

The proof of Theorem~\ref{BarStructure:thm:Existence} is now achieved.\qed

\medskip
Theorem~\ref{BarStructure:Existence:ModuleLevel} gives as a corollary:

\begin{mainthm}\label{BarStructure:Existence:FunctorLevel}
Suppose we have a morphism $\nabla_{\epsilon}: \QOp\rightarrow\End_{B_{\EOp}}$
so that $B_{\EOp}$ forms a $\QOp$-algebra in right $\EOp$-modules
as asserted in Theorem~\ref{BarStructure:Existence:ModuleLevel}.

Then
the bar complex $B(A) = \Sym_{\EOp}(B_{\EOp},A)$, $A\in{}_{\EOp}\E$,
becomes equipped with an induced $\QOp$-algebra structure
such that:
\begin{enumerate}
\item\label{ExistenceUniquenessFunctorLevel:Existence:Functoriality}
The operad $\QOp$ acts on $B(A)$ functorially in $A$.
\item\label{ExistenceUniquenessFunctorLevel:Existence:Restriction}
If $A$ is a commutative algebra,
then the action of $\QOp$ on $B(A)$ reduces to the standard action of the commutative operad on $B(A)$,
the action determined by the shuffle product of tensors.
\end{enumerate}
\end{mainthm}

\begin{proof}
To obtain assertion~(\ref{ExistenceUniquenessFunctorLevel:Existence:Functoriality}),
we use that the structure of a $\QOp$-algebra in right $\EOp$-modules
gives rise to a $\QOp$-algebra structure at the functor level.
Explicitly,
according to recollections of~\S\ref{Background:FunctorModules:FunctorsToAlgebras},
the evaluation morphism $\QOp(n)\otimes B_{\EOp}^{\otimes n}\rightarrow B_{\EOp}$
gives rise to an evaluation morphism at the functor level
\begin{equation*}
\QOp(n)\otimes\underbrace{\Sym_{\EOp}(B_{\EOp},A)}_{= B(A)}{}^{\otimes n}
\rightarrow\underbrace{\Sym_{\EOp}(B_{\EOp},A)}_{= B(A)}
\end{equation*}
so that the map $A\mapsto B(A)$ determines a functor from the category of $\EOp$-algebras
to the category of $\QOp$-algebras.

To obtain assertion~(\ref{ExistenceUniquenessFunctorLevel:Existence:Restriction}),
we use the relationship,
recalled in~\S\ref{Background:ExtensionRestriction:Functors},
between extensions and restrictions at the module and functor levels.
In the context of the theorem,
for a commutative algebra $A$,
we have a natural isomorphism
in the category of $\QOp$-algebras
\begin{equation*}
\underbrace{\Sym_{\COp}((B_{\EOp},\nabla_{\epsilon}),\epsilon^* A)}_{= (B(A),\nabla_{\epsilon})}
\simeq\Sym_{\COp}((B_{\EOp},\nabla_{\epsilon})\circ_{\EOp}\COp,A)
\end{equation*}
where the $\QOp$-algebra structure on the left-hand side comes from the bar module $B_{\EOp}$.
On the other hand,
we have a natural isomorphism
\begin{equation*}
\underbrace{\phi^*(\Sym_{\COp}((B_{\COp},\nabla_c),A))}_{= \phi^*(B(A),\nabla_c)}
\simeq\Sym_{\COp}(\phi^*(B_{\COp},\nabla_c),A),
\end{equation*}
where $\nabla_c$ represents the standard commutative algebra structure of the bar complex of $A$.
Hence,
if $(B_{\EOp},\nabla_{\epsilon})\circ_{\EOp}\COp\simeq\phi^*(B_{\COp},\nabla_c)$,
then we have a natural isomorphism of $\QOp$-algebras $(B(A),\nabla_{\epsilon})\simeq\phi^*(B(A),\nabla_c)$,
for all $A\in{}_{\COp}\E$.
\qed\end{proof}

\subsection{The uniqueness theorem}\label{BarStructure:Uniqueness}
In this subsection,
we prove that all solutions of the existence Theorem~\ref{BarStructure:thm:Existence}
yield equivalent objects in the homotopy category of $\QOp$-algebras in right $\EOp$-modules,
as well as equivalent structures on the bar construction at the functor level.

For this aim,
we use the morphism of right $\EOp$-modules
$\epsilon_{\sharp}: B_{\EOp}\rightarrow B_{\COp}$,
adjoint to the natural isomorphism
$\epsilon_{\flat}: B_{\EOp}\circ_{\EOp}\COp\xrightarrow{\simeq} B_{\COp}$
considered in Theorem~\ref{BarStructure:thm:Existence}.
By Proposition~\ref{BarConstruction:BarModule:BarModuleHomotopy}
this morphism $\epsilon_{\sharp}: B_{\EOp}\rightarrow B_{\COp}$
defines an acyclic fibration
since the augmentation of an $E_\infty$-operad $\epsilon: \EOp\rightarrow\COp$
forms itself an acyclic fibration in the category of operads.
Furthermore:

\begin{lemm}\label{BarStructure:Uniqueness:HomotopyType}
Suppose that the bar module $B_{\EOp}$
is equipped with the structure of a $\QOp$-algebra in right $\EOp$-modules
so that the natural isomorphism
$\epsilon_{\flat}: B_{\EOp}\circ_{\EOp}\COp\xrightarrow{\simeq} B_{\COp}$,
defines an isomorphism of $\QOp$-algebras in right $\COp$-modules
\begin{equation*}
\epsilon_{\flat}: (B_{\EOp},\nabla_{\epsilon})\circ_{\EOp}\COp\xrightarrow{\simeq}(B_{\COp},\nabla_c),
\end{equation*}
as asserted in Theorem~\ref{BarStructure:thm:Existence}.

Then
the morphism of right $\EOp$-modules
$\epsilon_{\sharp}: B_{\EOp}\rightarrow B_{\COp}$,
adjoint to
$\epsilon_{\flat}: B_{\EOp}\circ_{\EOp}\COp\xrightarrow{\simeq} B_{\COp}$,
defines a morphism of $\QOp$-algebras in right $\EOp$-modules
\begin{equation*}
\epsilon_{\sharp}: (B_{\EOp},\nabla_{\epsilon})\rightarrow(B_{\COp},\nabla_c),
\end{equation*}
and, hence, forms an acyclic fibration in that category.
\end{lemm}

\begin{proof}
In~\S\ref{Background:ExtensionRestriction:Algebras},
we recall that the extension and restriction functors $\psi_!: \M_{\ROp}\rightleftarrows\M_{\SOp} :\psi^*$
associated to any operad morphism $\psi: \ROp\rightarrow\SOp$
restrict to functors on $\POp$-algebras, for any operad $\POp$,
so that we have an adjunction relation:
\begin{equation*}
\psi_!: {}_{\POp}\M{}_{\ROp}\rightleftarrows{}_{\POp}\M{}_{\SOp} :\psi^*.
\end{equation*}
The lemma is an immediate corollary of this proposition.
\qed\end{proof}

This lemma gives immediately:

\begin{mainthm}\label{BarStructure:Uniqueness:ModuleLevel}
Suppose we have operad morphisms $\nabla_0,\nabla_1: \QOp\rightarrow\End_{B_{\EOp}}$
that provide the bar module $B_{\EOp}$ with the structure of a $\QOp$-algebra
in accordance with requirements (\ref{thm:Existence:ModuleAlgebra}-\ref{thm:Existence:CommutativeReduction})
of Theorem~\ref{BarStructure:thm:Existence}.

The algebras $(B_{\EOp},\nabla_0)$ and $(B_{\EOp},\nabla_1)$
are connected by weak-equivalences
\begin{equation*}
(B_{\EOp},\nabla_0)\xrightarrow{\sim}\,\cdot\,\xleftarrow{\sim}(B_{\EOp},\nabla_1)
\end{equation*}
in the category of $\QOp$-algebras in right $\EOp$-modules.\qed
\end{mainthm}

As usual in a model category,
the weak-equivalences
\begin{equation*}
(B_{\EOp},\nabla_0)\xrightarrow{\sim}\,\cdot\,\xleftarrow{\sim}(B_{\EOp},\nabla_1)
\end{equation*}
can be replaced by a chain of weak-equivalences of $\QOp$-algebras in right $\EOp$-modules
\begin{equation*}
(B_{\EOp},\nabla_0)\xleftarrow{\sim}\,\cdot\,\xrightarrow{\sim}\,\cdots\,\xrightarrow{\sim}(B_{\EOp},\nabla_1)
\end{equation*}
in which all intermediate objects are cofibrant as $\QOp$-algebras in right $\EOp$-modules,
and hence as right $\EOp$-modules
since any cofibrant algebra over a (cofibrant) operad $\QOp$
forms a cofibrant object in the underlying category (by~\cite[Corollary 5.5]{BergerMoerdijk}, \cite[Proposition 12.3.2]{FresseModules}).
Recall that the bar module $B_{\EOp}$
forms itself a cofibrant $\EOp$-module by Proposition~\ref{BarConstruction:BarModule:CofibrantStructure}.

In~\cite[\S 15]{FresseModules},
we prove that the natural transformation
\begin{equation*}
\Sym_{\EOp}(f,A): \Sym_{\EOp}(M,A)\xrightarrow{\sim}\Sym_{\EOp}(N,A)
\end{equation*}
induced by a weak-equivalence $f: M\xrightarrow{\sim} N$ such that $M,N$ are cofibrant right $\EOp$-modules
forms a weak-equivalence
for all $\EOp$-algebras $A$
which are cofibrant in the underlying category (see Theorem~15.1.A in \emph{loc. cit.}).
Accordingly, in our context, we obtain:

\begin{mainthm}\label{BarStructure:Uniqueness:FunctorLevel}
Suppose we have morphisms $\nabla_0,\nabla_1: \QOp\rightarrow\End_{B_{\EOp}}$,
as in Theorem~\ref{BarStructure:Existence:ModuleLevel},
that yield functorial $\QOp$-algebra structures on the bar construction $B(A)$
as in Theorem~\ref{BarStructure:Existence:FunctorLevel}.

The $\QOp$-algebras $(B(A),\nabla_0)$ and $(B(A),\nabla_1)$ can be connected by morphisms of $\QOp$-algebras
\begin{equation*}
(B(A),\nabla_0)\xleftarrow{\sim}\,\cdot\,\xrightarrow{\sim}\,\cdots\,\xrightarrow{\sim}(B(A),\nabla_1),
\end{equation*}
functorially in $A$,
and these morphisms are weak-equivalences whenever the $\EOp$-algebra~$A$
defines a cofibrant object in the underlying category $\E$.\qed
\end{mainthm}

\section{The categorical bar module}\label{CategoricalBarConstruction}

\subsection*{Introduction}
In the next section we prove that the bar construction $B(A)$,
equipped with the algebra structure of Theorem~\ref{BarStructure:Existence:FunctorLevel},
defines a model of the suspension in the homotopy category of $E_\infty$-algebras.
For this aim
we use a model of the suspension, defined in the general setting of pointed simplicial model categories
and yielded by a categorical version of the bar construction.

The purpose of this section
is to recall the definition of this categorical bar construction $C(A)$
in the context of algebras over an operad $\ROp$
and to define an $\ROp$-algebra in right $\ROp$-modules $C_{\ROp}$
such that $C(A) = \Sym_{\ROp}(C_{\ROp},A)$.
The plan of this section parallels the plan of~\S\ref{BarConstruction}
on the bar module $B_{\ROp}$.
In~\S\ref{CategoricalBarConstruction:CategoricalBarComplex},
we recall the definition of the categorical bar construction $C(A)$
in the context of algebras over operads,
where we take either the category of dg-modules $\E = \C$
or a category of right modules over an operad $\E = \M_{\SOp}$
as an underlying symmetric monoidal category;
in~\S\ref{CategoricalBarConstruction:OperadModuleBarConstruction},
we study the categorical bar construction of algebras in right modules over operads;
in~\S\ref{CategoricalBarConstruction:CategoricalBarModule},
we observe that the required $\ROp$-algebra in right $\ROp$-modules $C_{\ROp}$
is returned by the categorical bar construction of the $\ROp$-algebra in right $\ROp$-modules
formed by the operad itself.
Then we examine the functoriality of the construction $\ROp\mapsto C_{\ROp}$
and the homotopy invariance of the categorical bar module $C_{\ROp}$.

\subsection{Recollections: the categorical bar construction}\label{CategoricalBarConstruction:CategoricalBarComplex}
The categorical bar complex $C(A)$
is defined by the realization of a simplicial construction $\underline{C}(A)$
whose definition makes sense in any pointed category (explicitly, in any category equipped with a zero object $*$).
For our purpose,
we recall this definition in the context of algebras over a non-unitary operad $\POp$,
assumed to satisfy $\POp(0) = 0$,
and where the underlying category $\E$ is either the category of dg-modules itself $\E = \C$
or a category of right modules over an operad $\E = \M_{\SOp}$.
Note simply that the zero object of $\E$ is equipped with a $\POp$-algebra structure
if $\POp$ is a non-unitary operad
and defines obviously a zero object in ${}_{\POp}\E$.
Thus the category of $\POp$-algebras in $\E$,
where $\POp$ is any non-unitary operad,
is tautologically pointed.

Recall that $\Op_0$ denotes the category of non-unitary operads.

\subsubsection{The simplicial categorical bar complex}\label{CategoricalBarConstruction:CategoricalBarComplex:Simplicial}
To define the categorical bar construction $C_{\POp}(A)$ of an algebra $A\in{}_{\POp}\E$,
we form first a simplicial $\POp$-algebra $\underline{C}(A)$
such that
\begin{equation*}
\underline{C}(A)_n = A^{\vee n},
\end{equation*}
where $\vee$ denotes the categorical coproduct in the category of $\POp$-algebras in $\E$.
The faces and degeneracies of $\underline{C}(A)$
are defined explicitly by formulas
\begin{align*}
d_i & = \begin{cases} 0\vee A^{\vee n-1}, & \text{for $i = 0$}, \\
A^{\vee i-1}\vee\nabla\vee A^{\vee n-i-1}, & \text{for $i = 1,\dots,n-1$}, \\
A^{\vee n-1}\vee 0, & \text{for $i = n$},
\end{cases} \\
s_j & = A^{\vee j}\vee 0\vee A^{\vee n-j},\quad\text{for $j = 0,\dots,n$},
\end{align*}
where
$\nabla: A\vee A\rightarrow A$ denotes the codiagonal of $A$.

\subsubsection{On normalized complexes}\label{CategoricalBarConstruction:CategoricalBarComplex:NormalizedComplexes}
In the context of dg-modules $\E = \C$,
we use the standard normalized chain complex to associate a dg-module $N_*(\underline{C})$
to any simplicial dg-modules $\underline{C}$.

For a simplicial $\Sigma_*$-module $\underline{C}$,
the collection of normalized chain complexes $N_*(\underline{C}(n))$, $n\in\NN$,
defines a $\Sigma_*$-module $N_*(\underline{C})$ naturally associated to $\underline{C}$.
For a simplicial right $\ROp$-module $\underline{C}$,
we have an obvious isomorphism
\begin{equation*}
N_*(\underline{C})\circ\ROp\xrightarrow{\simeq} N_*(\underline{C}\circ\ROp),
\end{equation*}
so that $N_*(\underline{C})$
inherits the structure of a right $\ROp$-module
and defines an object of~$\M_{\ROp}$.

In our constructions,
we use the classical Eilenberg-Mac Lane equivalence,
which gives a natural morphism
\begin{equation*}
N_*(\underline{C})\otimes N_*(\underline{D})\xrightarrow{\EM} N_*(\underline{C}\otimes\underline{D}),
\end{equation*}
for all simplicial dg-modules $\underline{C},\underline{D}$.
In the context of a category of right modules over an operad $\E = \M_{\ROp}$,
we have termwise Eilenberg-Mac Lane morphisms
\begin{equation*}
\Sigma_{r}\otimes_{\Sigma_s\times\Sigma_t} N_*(\underline{C}(s))\otimes N_*(\underline{D}(t))
\xrightarrow{\EM} N_*(\Sigma_{r}\otimes_{\Sigma_s\times\Sigma_t}\underline{C}(s)\otimes\underline{D}(t)),
\end{equation*}
inherited from dg-modules,
which assemble to give an Eilenberg-Mac Lane morphism in~$\M_{\ROp}$
\begin{equation*}
N_*(\underline{C})\otimes N_*(\underline{D})\xrightarrow{\EM} N_*(\underline{C}\otimes\underline{D}),
\end{equation*}
and similarly as regards the external tensor product in~$\M_{\ROp}$.

In all cases,
if $\underline{C}\equiv C$ is a constant simplicial object,
then the Eilenberg-Mac Lane morphism is identified with a natural isomorphism
\begin{equation*}
C\otimes N_*(\underline{D})\simeq N_*(C\otimes\underline{D}).
\end{equation*}

\subsubsection{The normalized categorical bar construction}\label{CategoricalBarConstruction:CategoricalBarComplex:Realization}
The categorical bar construction $C(A)$
is defined by the normalized chain complex
\begin{equation*}
C(A) = N_*(\underline{C}(A)).
\end{equation*}

This object is equipped with the structure of a $\POp$-algebra,
like the normalized chain complex of any simplicial algebra over an operad.
Formally,
we have evaluation products
\begin{equation*}
\POp(n)\otimes C(A)^{\otimes n}\rightarrow C(A)
\end{equation*}
defined by the composite of the Eilenberg-Mac Lane equivalences
\begin{equation*}
N_*(\underline{C}(A))^{\otimes n}\xrightarrow{\EM} N_*(\underline{C}(A)^{\otimes n})
\end{equation*}
with the morphisms
\begin{equation*}
\POp(n)\otimes N_*(\underline{C}(A)^{\otimes n}) = N_*(\POp(n)\otimes\underline{C}(A)^{\otimes n})
\rightarrow N_*(\underline{C}(A))
\end{equation*}
induced by the evaluation product of~$\underline{C}(A)$.

\subsection{The categorical bar construction of algebras in right modules over operads}\label{CategoricalBarConstruction:OperadModuleBarConstruction}
In this section,
we study the categorical bar construction of $\POp$-algebras
in right modules over an operad $\ROp$.
In this context,
the categorical bar construction $N\mapsto C(N)$ returns a $\POp$-algebra in right $\ROp$-modules.
As in~\S\ref{BarConstruction:AinfinityBarComplex:FunctorBarModule},
we determine the functor $\Sym_{\ROp}(C(N)): {}_{\ROp}\E\rightarrow{}_{\POp}\E$
associated to this object $C(N)\in{}_{\POp}\M{}_{\ROp}$.

In the context of the standard bar construction,
we use that the functor $\Sym_{\ROp}: \M{}_{\ROp}\rightarrow\Func{}_{\ROp}$ preserves tensor products
to identify the functor $A\mapsto \Sym_{\ROp}(B(N),A)$
associated to the bar complex of a $\KOp$-algebra in right $\ROp$-modules
with the bar complex $B(\Sym_{\ROp}(N,A))$
of the $\KOp$-algebra $\Sym_{\ROp}(N,A)\in{}_{\KOp}\E$.
Similarly,
as the functor $\Sym_{\ROp}: {}_{\POp}\M{}_{\ROp}\rightarrow{}_{\POp}\Func{}_{\ROp}$
preserves colimits of $\POp$-algebras (see~\S\ref{Background:FunctorModules:FunctorsToAlgebras}),
we obtain:

\begin{lemm}\label{CategoricalBarConstruction:OperadModuleBarConstruction:SimplicialFunctor}
Let ${}_{\POp}\E^{\Delta}$ be the category of simplicial $\POp$-algebras.
Let $N$ be a $\POp$-algebra in right $\ROp$-modules.

The functor $\Sym_{\ROp}(\underline{C}(N)): {}_{\ROp}\E\rightarrow{}_{\POp}\E^{\Delta}$
associated to the simplicial categorical bar construction of $N$
satisfies the identity
\begin{equation*}
\Sym_{\ROp}(\underline{C}(N),A) = \underline{C}(\Sym_{\ROp}(N,A)),
\end{equation*}
for all $A\in{}_{\ROp}\E$,
where on the right-hand side we consider the simplicial categorical bar complex of the $\POp$-algebra $\Sym_{\ROp}(N,A)$
associated to $A\in{}_{\ROp}\E$ by the functor $\Sym_{\ROp}(N): {}_{\ROp}\E\rightarrow {}_{\POp}\E$
defined by $N$.\qed
\end{lemm}

As the normalized chain complex $N_*(\underline{C})$
of a simplicial object $\underline{C}$
is defined by a cokernel
and the functor $M\mapsto\Sym_{\ROp}(M)$ preserve colimits in right $\ROp$-modules,
we have a natural isomorphism $\Sym_{\ROp}(N_*(\underline{M}),A)\simeq N_*(\Sym_{\ROp}(\underline{M},A))$,
for all $A\in {}_{\ROp}\E$.
This isomorphism commutes with Eilenberg-Mac Lane equivalences
in the sense that the coherence diagram
\begin{equation*}
\xymatrix@!C=3cm{ & \Sym_{\ROp}(N_*(\underline{C}),A)\otimes \Sym_{\ROp}(N_*(\underline{D}),A)\ar[dr]^{\simeq} & \\
\Sym_{\ROp}(N_*(\underline{C})\otimes N_*(\underline{D}),A)\ar[ur]^{\simeq}\ar[d]_{\EM} & &
N_*(\Sym_{\ROp}(\underline{C},A))\otimes N_*(\Sym_{\ROp}(\underline{D},A))\ar[d]_{\EM} \\
\Sym_{\ROp}(N_*(\underline{C}\otimes\underline{D}),A)\ar[dr]_{\simeq} & & N_*(\Sym_{\ROp}(\underline{C},A)\otimes\Sym(\underline{D},A)) \\
& N_*(\Sym_{\ROp}(\underline{C}\otimes\underline{D},A))\ar[ur]_{\simeq} & }
\end{equation*}
commutes.
As a consequence,
if $\underline{C}$ is a simplicial $\POp$-algebra in right $\ROp$-modules,
then the functor identity $\Sym_{\ROp}(N_*(\underline{C}),A) = N_*(\Sym_{\ROp}(\underline{C},A))$
holds in the category of $\POp$-algebras.

From these observations,
we conclude:

\begin{prop}\label{CategoricalBarConstruction:OperadModuleBarConstruction:Functor}
Let $N$ be a $\POp$-algebra in right $\ROp$-modules.

The functor $\Sym_{\ROp}(C(N)): {}_{\ROp}\E\rightarrow{}_{\POp}\E$
associated to the categorical bar construction of $N$
satisfies the relation
\begin{equation*}
\Sym_{\ROp}(C(N),A)\simeq C(\Sym_{\ROp}(N,A)),
\end{equation*}
for all $A\in{}_{\ROp}\E$,
where on the right-hand side we consider the categorical bar complex of the $\POp$-algebra $\Sym_{\ROp}(N,A)$
associated to $A\in{}_{\ROp}\E$ by the functor $\Sym_{\ROp}(N): {}_{\ROp}\E\rightarrow {}_{\POp}\E$
defined by $N$.\qed
\end{prop}

\begin{remark}
In~\S\ref{BarConstruction:OperadModuleBarComplex},
we observe that the functor $N\mapsto B(N)$ commutes with extensions and restrictions of structure on the right.
The same assertion holds for the functor $N\mapsto C(N)$
defined by the categorical bar construction
just because both functors $\psi_!: {}_{\POp}\M{}_{\ROp}\rightleftarrows{}_{\POp}\M{}_{\SOp} :\psi^*$
preserve coproducts.
The functor $N\mapsto C(N)$
also commutes with extensions of structure on the left,
but not with restrictions of structure on the left
since this latter operation does not preserve coproducts.
Nevertheless,
we still have a natural morphism $C(\phi^* N)\rightarrow\phi^* C(N)$
induced by the natural transformations $(\phi^* N)^{\vee n}\rightarrow\phi^*(N^{\vee n})$.
\end{remark}

\renewcommand{\themainthm}{\thesubsection.\Alph{mainthm}}

\subsection{The categorical bar module}\label{CategoricalBarConstruction:CategoricalBarModule}
The categorical bar module of an operad $\ROp$, like the bar module of~\S\ref{BarConstruction:BarModule},
is the categorical bar construction
of the $\ROp$-algebra in right $\ROp$-modules formed by the operad itself.
For the sake of coherence,
we use the notation $C_{\ROp}$ for this categorical bar module $C_{\ROp} = C(\ROp)$
and
we set similarly $\underline{C}{}_{\ROp} = \underline{C}(\ROp)$.

In~\S\ref{Background:FunctorModules:FunctorsToAlgebras},
we recall that $\Sym_{\ROp}(\ROp): {}_{\ROp}\E\rightarrow {}_{\ROp}\E$
represents the identity functor on the category of $\ROp$-algebras.
Hence,
Proposition~\ref{CategoricalBarConstruction:OperadModuleBarConstruction:Functor}
gives:

\begin{mainprop}\label{CategoricalBarConstruction:CategoricalBarModule:FunctorCategoricalBarModule}
The functor $\Sym_{\ROp}(C_{\ROp}): {}_{\ROp}\E\rightarrow{}_{\ROp}\E$ associated to $C_{\ROp}$
is naturally isomorphic to the categorical bar construction $A\mapsto C(A)$
in the category of $\ROp$-algebras.\qed
\end{mainprop}

As in~\S\ref{BarConstruction:BarModule}, we examine the structure of $C_{\ROp}$
and the functoriality of the construction $\ROp\mapsto C_{\ROp}$.

The categorical bar module $C_{\ROp}$ does not form a cofibrant object in right $\ROp$-modules,
unlike the bar module $B_{\ROp}$,
but we prove that $C_{\ROp}$ is cofibrant as a $\Sigma_*$-module
provided that the operad $\ROp$ is so
(according to our usual convention, we say that $C_{\ROp}$ is $\Sigma_*$-cofibrant).
Thus,
we forget right module structures
and
we examine the $\ROp$-algebra in $\Sigma_*$-modules
underlying the categorical bar module $C_{\ROp}$.
For the simplicial categorical bar module $\underline{C}{}_{\ROp}$,
we obtain:

\begin{lemm}\label{CategoricalBarConstruction:CategoricalBarModule:UnderlyingLeftModule}
We have an identity
\begin{equation*}
\underline{C}{}_{\ROp} = \ROp(\underline{C}(\IOp)),
\end{equation*}
where $\ROp(\underline{C}(\IOp))$ represents the free $\ROp$-algebra
on the categorical bar construction of the unit $\Sigma_*$-module $\IOp$
in the category of $\Sigma_*$-modules.
\end{lemm}

\begin{proof}
By construction,
the forgetful functor $U: \M{}_{\ROp}\rightarrow\M$ preserves enriched monoidal category structures.
By~\cite[Proposition 3.3.3]{FresseModules},
this assertion implies that the forgetful functor $U: {}_{\ROp}\M{}_{\ROp}\rightarrow{}_{\ROp}\M$,
from the category of $\ROp$-algebras in right $\ROp$-modules
to the category of $\ROp$-algebras in $\Sigma_*$-modules,
preserves colimits.
As a consequence,
we obtain that $\underline{C}_{\ROp}$
agrees with the categorical bar construction of $\ROp$
in $\Sigma_*$-modules.

Observe that the operad $\ROp$ forms a free object in the category of $\ROp$-algebras in $\Sigma_*$-modules:
we have explicitly $\ROp = \ROp\circ\IOp = \ROp(\IOp)$.
By adjunction,
a coproduct of free objects satisfies the relation $\ROp(M)\vee\ROp(N) = \ROp(M\oplus N)$,
for all $M,N\in\M$.
Hence,
we obtain readily:
\begin{equation*}
(\underline{C}{}_{\ROp})_n = \ROp^{\vee n} = \ROp(\IOp)^{\vee n} = \ROp(\IOp^{\oplus n}),
\end{equation*}
for all $n\in\NN$.
The determination of faces and degeneracies of the categorical bar construction
is also formal from the universal property of free objects,
so that we obtain the conclusion of the lemma.
\qed\end{proof}

As a byproduct,
we obtain:

\begin{prop}\label{CategoricalBarConstruction:CategoricalBarModule:UnderlyingSigmaModule}
The categorical bar module $C_{\ROp}$ is $\Sigma_*$-cofibrant
if the operad $\ROp$ is so.
\end{prop}

\begin{proof}
The assumption about the operad $\ROp$ implies that simplicial $\ROp$-algebras
form a model category (see references of~\S\ref{Background:OperadAlgebras:AlgebraModelCategories}).
Lemma~\ref{CategoricalBarConstruction:CategoricalBarModule:UnderlyingLeftModule}
implies that the simplicial categorical bar module $\underline{C}_{\ROp}$
forms a cofibrant simplicial $\ROp$-algebra in $\Sigma_*$-modules.
By~\cite[Corollary 5.5]{BergerMoerdijk},~\cite[Proposition 12.3.2]{FresseModules},
this assertion implies
that $\underline{C}_{\ROp}$ is cofibrant in the underlying category of simplicial $\Sigma_*$-modules,
and hence that the normalized chain complex $C_{\ROp} = N_*(\underline{C}_{\ROp})$
associated to $\underline{C}_{\ROp}$ is cofibrant as a $\Sigma_*$-module.
\qed\end{proof}

\begin{construction}[Functoriality of the categorical bar module]\label{CategoricalBarConstruction:CategoricalBarModule:Functoriality}
In~\S\ref{BarConstruction:BarModule},
we observe that a morphism of operads under $\KOp$ gives rise to a morphism of right $\ROp$-modules
$\psi_{\sharp}: B_{\ROp}\rightarrow B_{\SOp}$.
In this paragraph,
we check that a morphism of operads gives rise to an analogous morphism of $\ROp$-algebras in right $\ROp$-modules
\begin{equation*}
\psi_{\sharp}: C_{\ROp}\rightarrow\psi^* C_{\SOp},
\end{equation*}
where $\psi^* C_{\SOp}$ refers to the $\ROp$-algebra in right $\ROp$-modules
obtained by a two-sided restriction of $C_{\SOp}\in{}_{\SOp}\M{}_{\SOp}$.
We prove next that $\psi_{\sharp}$ defines a weak-equivalence (respectively, a fibration) if $\psi$ is so.

Formally,
we use that $\psi$ determines a two-sided restriction functor $\psi^*: {}_{\SOp}\M{}_{\SOp}\rightarrow{}_{\ROp}\M{}_{\ROp}$
and
the operad morphism $\psi: \ROp\rightarrow\SOp$ defines a morphism $\psi: \ROp\rightarrow\psi^*\SOp$
in the category of $\ROp$-algebras in right $\ROp$-modules.
As a consequence,
by functoriality of the categorical bar construction $N\mapsto C(N)$,
we obtain that $\psi: \ROp\rightarrow\SOp$
induces a natural morphism of $\ROp$-algebras in right $\ROp$-modules
\begin{equation*}
C(\ROp)\xrightarrow{C(\psi)} C(\psi^*\SOp).
\end{equation*}

On the other hand,
for any algebra $N\in {}_{\SOp}\M{}_{\SOp}$,
we have a morphism
$\underline{C}(\psi^* N)\rightarrow\psi^*\underline{C}(N)$
induced by the natural transformation $(\psi^* N)^{\vee n}\rightarrow\psi^*(N^{\vee n})$.
As a consequence,
we have a natural morphism
\begin{equation*}
C(\psi^* N)\xrightarrow{\psi_{\sharp}}\psi^* C(N).
\end{equation*}
between the categorical bar complex of $N\in{}_{\SOp}\M{}_{\SOp}$ and the categorical bar complex of $\psi^* N\in{}_{\ROp}\M{}_{\ROp}$.
Our morphism $\psi_{\sharp}: C_{\ROp}\rightarrow\psi^* C_{\SOp}$
is given by the composite:
\begin{equation*}
C(\ROp)\xrightarrow{C(\psi)} C(\psi^*\SOp)\xrightarrow{\psi_{\sharp}}\psi^* C(\SOp).
\end{equation*}
\end{construction}

If we forget right module structures,
then we obtain readily:

\begin{obsv}\label{CategoricalBarConstruction:CategoricalBarModule:FunctorialityRepresentation}
The morphism
\begin{equation*}
\underline{C}{}_{\ROp}\xrightarrow{\psi_{\sharp}}\psi^*\underline{C}{}_{\SOp}
\end{equation*}
associated to an operad morphism $\psi: \ROp\rightarrow\SOp$
is given dimensionwise by the natural morphism of free objects
\begin{equation*}
\ROp(\underline{C}(\IOp))\xrightarrow{\psi(\underline{C}(\IOp))}\SOp(\underline{C}(\IOp))
\end{equation*}
induced by $\psi: \ROp\rightarrow\SOp$.
\end{obsv}

We use this observation
to prove:

\begin{lemm}\label{CategoricalBarConstruction:CategoricalBarModule:Fibrations}
If $\psi: \ROp\rightarrow\SOp$ is a weak-equivalence (respectively a fibration) of operads,
then the morphism
\begin{equation*}
\psi_{\sharp}: C_{\ROp}\rightarrow\psi^* C_{\SOp}
\end{equation*}
defines a weak-equivalence (respectively a fibration) in ${}_{\ROp}\M{}_{\ROp}$, for all non-unitary operads $\ROp,\SOp\in\Op_0$.
\end{lemm}

\begin{proof}
Since all forgetful functors create weak equivalences and fibrations,
we can forget right module structures
in the proof of this lemma
and
we can use the representation of Observation~\ref{CategoricalBarConstruction:CategoricalBarModule:FunctorialityRepresentation}.

We deduce immediately
from the form of the free $\ROp$-algebra
\begin{equation*}
\ROp(\underline{C}(\IOp)) = \Sym(\ROp,\underline{C}(\IOp)) = \ROp\circ\underline{C}(\IOp)
\end{equation*}
that the morphism of simplicial dg-modules
$\psi(\underline{C}(\IOp)) = \psi\circ\underline{C}(\IOp): \ROp\circ\underline{C}(\IOp)\rightarrow\SOp\circ\underline{C}(\IOp)$
induced by a surjective morphism of dg-operads $\phi: \ROp\rightarrow\SOp$
is surjective as well.
As a byproduct,
so is the morphism induced by $\psi(\underline{C}(\IOp))$
on normalized chain complexes.
Thus
we conclude that the morphism
$\psi_{\sharp}: C_{\ROp}\rightarrow\psi^* C_{\SOp}$
induced by a fibration of dg-operads
forms a fibration as well.

Recall that the composition product of $\Sigma_*$-modules $M\circ N$
preserves weak-equivalences in $M$,
provided that $N(0) = 0$ and the modules $N(r)$, $r>0$,
are cofibrant in dg-modules
(see~\cite[\S 2.3]{FressePartitions}, see also~\cite[\S 11.6]{FresseModules}).
From this assertion, we deduce that the morphism of simplicial $\Sigma_*$-modules
$\psi(\underline{C}(\IOp)) = \psi\circ\underline{C}(\IOp): \ROp\circ\underline{C}(\IOp)\rightarrow\SOp\circ\underline{C}(\IOp)$
induced by a weak-equivalence of dg-operads $\psi: \ROp\rightarrow\SOp$
defines a weak-equivalence,
and so does the morphism induced by $\psi(\underline{C}(\IOp))$
on normalized chain complexes.
Hence,
we conclude that the morphism
$\psi_{\sharp}: C_{\ROp}\rightarrow\psi^* C_{\SOp}$
induced by a weak-equivalence of dg-operads
$\psi: \ROp\xrightarrow{\sim}\SOp$
forms a weak-equivalence.
\qed\end{proof}

\renewcommand{\themainthm}{\thesection.\Alph{mainthm}}

\section{The homotopy interpretation of the bar construction}\label{HomotopyInterpretation}

\subsection*{Introduction}
In this section,
we prove that, for cofibrant algebras over $E_\infty$-operads,
the usual bar construction $B(A)$, equipped with the algebra structure given by Theorem~\ref{BarStructure:Existence:FunctorLevel},
is equivalent to the categorical bar construction $C(A)$
as an $E_\infty$-algebra.
Then we use that the categorical bar construction $C(A)$
is equivalent to the suspension $\Sigma A$
in the homotopy categories of algebras over an operad
to conclude:

\begin{mainthm}\label{HomotopyInterpretation:thm:HomotopyInterpretation}
Suppose that $\EOp$ forms itself a cofibrant $E_\infty$-operad
and set $\QOp = \EOp$.

Assume
that the bar complex $B(A)$ is equipped with the structure of an $\EOp$-algebra, for all $A\in{}_{\EOp}\E$,
and
that this structure is realized at the module level,
as stated in Theorem~\ref{BarStructure:Existence:FunctorLevel}.
Then
we have natural $\EOp$-algebra equivalences
\begin{equation*}
B(A)\xleftarrow{\sim}\,\cdot\,\xrightarrow{\sim}\,\cdots\,\xrightarrow{\sim}\Sigma A
\end{equation*}
that connect $B(A)$ to the suspension of $A$ in the model category of $\EOp$-algebras,
for all cofibrant $\EOp$-algebras $A$.\qed
\end{mainthm}

This theorem can easily be generalized to include the case where the operad $\EOp$ is not itself cofibrant (see~\S\ref{HomotopyInterpretation:SuspensionModel}).

\medskip
Again we realize the equivalence between $B(A)$ and $C(A)$ at the module level.
To be explicit,
let $\EOp$ be any $E_\infty$-operad (possibly not cofibrant), let $\QOp$ be any cofibrant $E_\infty$-operad,
and assume that the bar module $B_{\EOp}$ is equipped with the structure of a $\QOp$-algebra in right $\EOp$-modules,
as asserted in Theorem~\ref{BarStructure:thm:Existence}.
Recall that the categorical bar module $B_{\EOp}$ forms an $\EOp$-algebra in right $\EOp$-modules.
Since $\QOp$ is supposed to be cofibrant,
we can pick an operad morphism
in the lifting diagram
\begin{equation*}
\xymatrix{ & \EOp\ar@{->>}[d]_{\sim}^{\epsilon} \\ \QOp\ar@{.>}[ur]^{\psi}\ar[r]_{\phi} & \COp }
\end{equation*}
to make any $\EOp$-algebra in right $\EOp$-modules
into a $\QOp$-algebra in right $\EOp$-modules
by restriction of structure.
In~\S\ref{HomotopyInterpretation:EinfinityBarConstructionEquivalence},
we check that $B_{\EOp}$ and $C_{\EOp}$ define equivalent objects
in the homotopy category of $\QOp$-algebras in right $\EOp$-modules.
Thus we have a chain of weak-equivalences of $\QOp$-algebras in right $\EOp$-modules
\begin{equation*}
B_{\EOp}\xleftarrow{\sim}\,\cdot\,\xrightarrow{\sim}\,\cdots\,\xrightarrow{\sim} C_{\EOp}.
\end{equation*}

In~\S\ref{HomotopyInterpretation:SuspensionModel},
we use a theorem of~\cite[\S 15]{FresseModules}
to obtain that these weak-equivalences give rise to weak-equivalences at the functor level
\begin{equation*}
\underbrace{\Sym_{\EOp}(B_{\EOp},A)}_{=B(A)}\xleftarrow{\sim}\,\cdot\,\xrightarrow{\sim}\,\cdots\,\xrightarrow{\sim} \underbrace{\Sym_{\EOp}(C_{\EOp},A)}_{=C(A)}
\end{equation*}
for all cofibrant $\EOp$-algebras $A\in{}_{\EOp}\E$
and our conclusion follows.

\subsection{The equivalence of bar constructions}\label{HomotopyInterpretation:EinfinityBarConstructionEquivalence}
First we prove the existence of equivalences between the bar modules $B_{\EOp}$ and $C_{\EOp}$
associated to an $E_\infty$-operad $\EOp$.
This result is a consequence of the following observation:

\begin{lemm}\label{HomotopyInterpretation:EinfinityBarConstructionEquivalence:CommutativeCaseIdentity}
For the commutative operad $\COp$,
we have an identity of $\COp$-algebras in right $\COp$-modules $B_{\COp} = C_{\COp}$.
\end{lemm}

\begin{proof}
This observation is a consequence of the definition of the coproduct in the category of non-unitary commutative algebras.
Explicitly, for non-unitary commutative algebras in dg-modules,
and more generally in any symmetric monoidal category,
we have an identity:
$A\vee B = A\oplus B\oplus A\otimes B$.
As a consequence,
for the simplicial categorical bar complex $\underline{C}(N)$ of any commutative algebra $N$ in right $\COp$-modules,
we obtain
\begin{equation*}
\underline{C}(N)_n = (N\otimes\dots\otimes N)\oplus(\text{degeneracies}).
\end{equation*}
Thus,
at the level of normalized chain complexes,
we obtain the relation
$C(N) = N_*(\underline{C}(N)) = B(N)$.
The case $N = \COp$ gives the announced identity $B_{\COp} = C_{\COp}$.
\qed\end{proof}

Roughly, for an $E_\infty$-operad $\EOp$, we lift the isomorphism of this lemma to a weak-equivalence
of $E_\infty$-algebras in right $\EOp$-modules.

Suppose we have a cofibrant $E_\infty$-operad $\QOp$
together with an augmentation $\phi: \QOp\xrightarrow{\sim}\COp$.
Assume that the bar module $B_{\EOp}$
is equipped with the structure of a $\QOp$-algebra in right $\EOp$-modules
as asserted in Theorem~\ref{BarStructure:thm:Existence}.
In Proposition~\ref{BarStructure:Uniqueness:HomotopyType},
we observe that the obtained $\QOp$-algebra $B_{\EOp}$ is endowed with a weak-equivalence
\begin{equation*}
\epsilon_{\sharp}: B_{\EOp}\xrightarrow{\sim}\phi^* B_{\COp}
\end{equation*}
in the category of $\QOp$-algebras in right $\EOp$-modules.

On the other hand,
we observe in Proposition~\ref{CategoricalBarConstruction:CategoricalBarModule:Fibrations},
that the categorical bar module $C_{\EOp}$ is endowed with a weak-equivalence
\begin{equation*}
\epsilon_{\sharp}: C_{\EOp}\xrightarrow{\sim}\epsilon^* C_{\COp}
\end{equation*}
in the category of $\EOp$-algebras in right $\EOp$-modules.
As explained in the introduction of this section,
since $\QOp$ is cofibrant,
we can pick a lifting
in the operad diagram
\begin{equation*}
\xymatrix{ & \EOp\ar@{->>}[d]^{\sim} \\ \QOp\ar[r]_{\phi}\ar@{.>}[ur]^{\psi} & \COp }
\end{equation*}
to obtain a morphism $\psi: \QOp\rightarrow\EOp$
in $\Op_0/\COp$.
By restriction of structure,
the equivalence $\epsilon_{\sharp}: C_{\EOp}\xrightarrow{\sim}\epsilon^* C_{\COp}$
gives rise to an equivalence
\begin{equation*}
\psi^*(\epsilon_{\sharp}): \psi^* C_{\EOp}\xrightarrow{\sim}\psi^*\epsilon^* C_{\COp} = \phi^* C_{\COp}
\end{equation*}
in the category of $\QOp$-algebras in right $\EOp$-modules.

Therefore, we obtain:

\begin{thm}\label{HomotopyInterpretation:EinfinityBarConstructionEquivalence:ModuleLevel}
Assume that the bar module $B_{\EOp}$
is equipped with the structure of a $\QOp$-algebra in right $\EOp$-modules
as in Theorem~\ref{BarStructure:thm:Existence}.
Then we have weak-equivalences
\begin{equation*}
B_{\EOp}\xrightarrow{\sim}\phi^* B_{\COp} = \phi^* C_{\COp}\xleftarrow{\sim}\psi^* C_{\EOp}
\end{equation*}
in the category of $\QOp$-algebras in right $\EOp$-modules.\qed
\end{thm}

Again,
we can use model category structures
to replace the weak-equivalences
\begin{equation*}
B_{\EOp}\xrightarrow{\sim}\,\cdot\,\xleftarrow{\sim}\psi^* C_{\EOp}
\end{equation*}
by a chain of weak-equivalences
\begin{equation*}
B_{\EOp}\xleftarrow{\sim}\,\cdot\,\xrightarrow{\sim}\,\cdots\,\xrightarrow{\sim}\psi^* C_{\EOp}
\end{equation*}
in which all intermediate objects are cofibrant objects of the category of $\QOp$-algebras in right $\EOp$-modules.
Recall that a right $\EOp$-module $M$ is called $\Sigma_*$-cofibrant, like an operad,
if $M$ is cofibrant as a $\Sigma_*$-module.
By~\cite[Proposition 14.1.1]{FresseModules},
any cofibrant right $\EOp$-module is $\Sigma_*$-cofibrant
since the $E_\infty$-operad $\EOp$ is supposed to be $\Sigma_*$-cofibrant.
Accordingly, the bar module $B_{\EOp}$ is $\Sigma_*$-cofibrant.
The categorical bar module $C_{\EOp}$ is also $\Sigma_*$-cofibrant
by Proposition~\ref{CategoricalBarConstruction:CategoricalBarModule:UnderlyingSigmaModule}.
Since a cofibrant $\QOp$-algebra in right $\EOp$-modules
forms a cofibrant object in the underlying category of right $\EOp$-modules
by~\cite[Corollary 5.5]{BergerMoerdijk}, \cite[Proposition 12.3.2]{FresseModules},
and hence a $\Sigma_*$-cofibrant module by~\cite[Proposition 14.1.1]{FresseModules},
we conclude that all objects in our chain of weak-equivalences are $\Sigma_*$-cofibrant.
At the functor level,
we obtain that these weak-equivalences give rise to:

\begin{thm}\label{HomotopyInterpretation:EinfinityBarConstructionEquivalence:FunctorLevel}
The bar construction $B(A)$
is connected to the categorical bar construction $C(A)$
by natural weak-equivalences of $\QOp$-algebras
\begin{equation*}
B(A)\xleftarrow{\sim}\,\cdot\,\xrightarrow{\sim}\,\cdots\,\xrightarrow{\sim}\psi^* C(A),
\end{equation*}
for all cofibrant $\EOp$-algebras $A$,
where we use a restriction of structure to make the $\EOp$-algebra $C(A)$
into a $\QOp$-algebra.\qed
\end{thm}

\begin{proof}
In~\cite[\S 15]{FresseModules},
we prove that a weak-equivalence $\phi: M\xrightarrow{\sim} N$
between $\Sigma_*$-cofibrant right $\ROp$-modules $M$ and $N$
induces a weak-equivalence at the functor level:
\begin{equation*}
\Sym_{\ROp}(\phi,A): \Sym_{\ROp}(M,A)\xrightarrow{\sim} \Sym_{\ROp}(N,A),
\end{equation*}
for all cofibrant $\ROp$-algebras $A$.
Accordingly,
the morphisms
\begin{equation*}
B_{\EOp}\xleftarrow{\sim}\,\cdot\,\xrightarrow{\sim}\,\cdots\,\xrightarrow{\sim}\psi^* C_{\EOp}
\end{equation*}
induce weak-equivalences of $\QOp$-algebras
\begin{equation*}
\Sym_{\EOp}(B_{\EOp},A)\xleftarrow{\sim}\,\cdot\,\xrightarrow{\sim}\,\cdots\,\xrightarrow{\sim}\Sym_{\EOp}(\psi^* C_{\EOp},A).
\end{equation*}
for all cofibrant $\EOp$-algebras $A$.

Recall that the functor $N\mapsto\Sym_{\EOp}(N)$ commutes with restrictions of structure on the left.
Therefore we have weak-equivalences between
\begin{equation*}
B(A) = \Sym_{\EOp}(B_{\EOp},A)
\end{equation*}
and
\begin{equation*}
\psi^* C(A) = \psi^*\Sym_{\EOp}(C_{\EOp},A) = \Sym_{\EOp}(\psi^* C_{\EOp},A)
\end{equation*}
as required.
\qed\end{proof}

\subsection{The equivalence with suspensions}\label{HomotopyInterpretation:SuspensionModel}
The next assertion is proved in~\cite{Mandell}
(in the context of dg-modules, but the generalization to any category over dg-modules $\E$ is straightforward):

\begin{fact}[{See~\cite[\S 3, \S 14]{Mandell}}]\label{HomotopyInterpretation:SuspensionModel:EquivalenceToSuspension}
Assume that $\POp$ is a $\Sigma_*$-cofibrant operad in dg-modules
so that the category of $\POp$-algebras in $\E$
forms a semi-model category.

For every cofibrant $\POp$-algebra in $\E$,
the $\POp$-algebra $C(A)$ is connected to $\Sigma A$, the suspension of $A$ in the model category of $\POp$-algebras in $\E$,
by weak-equivalences of $\POp$-algebras
\begin{equation*}
C(A)\xleftarrow{\sim}\,\cdot\,\xrightarrow{\sim}\,\cdots\,\xrightarrow{\sim}\Sigma A,
\end{equation*}
functorially in $A$.
\end{fact}

From this assertion and Theorem~\ref{HomotopyInterpretation:EinfinityBarConstructionEquivalence:FunctorLevel} we conclude:

\begin{mainthm}[Claim of Theorem~\ref{HomotopyInterpretation:thm:HomotopyInterpretation}]
Suppose that $\EOp$ is a cofibrant $E_\infty$-operad
and set $\QOp = \EOp$.

Assume
that the bar complex $B(A)$ is equipped with the structure of an $\EOp$-algebra, for all $A\in{}_{\EOp}\E$,
and
that this structure is realized at the module level,
as stated in Theorem~\ref{BarStructure:thm:Existence}.
Then
we have natural $\EOp$-algebra equivalences
\begin{equation*}
B(A)\xleftarrow{\sim}\,\cdot\,\xrightarrow{\sim}\,\cdots\,\xrightarrow{\sim}\Sigma A
\end{equation*}
that connect $B(A)$ to the suspension of $A$ in the model category of $\EOp$-algebras,
for all cofibrant $\EOp$-algebras $A$.\qed
\end{mainthm}

To complete this result,
recall that the extension and restriction functors
\begin{equation*}
\phi_!: {}_{\POp}\E\rightleftarrows{}_{\QOp}\E :\phi^*.
\end{equation*}
associated to a weak-equivalence of $\Sigma_*$-cofibrant operads $\phi: \POp\rightarrow\QOp$
define Quillen adjoint equivalences of model categories (see~\cite{BergerMoerdijk} or~\cite[\S 16]{FresseModules}).
As a byproduct,
Theorem~\ref{HomotopyInterpretation:thm:HomotopyInterpretation} can be generalized to cover the case
where the $E_\infty$-operad $\EOp$ is not cofibrant as an operad.
In this context,
we obtain weak-equivalences of $\QOp$-algebras
\begin{equation*}
B(A)\xleftarrow{\sim}\,\cdot\,\xrightarrow{\sim}\,\cdots\,\xrightarrow{\sim}\psi^*(\Sigma A),
\end{equation*}
for all cofibrant $\EOp$-algebras $A$,
where $\Sigma A$ is the suspension of $A$ in the model category of $\EOp$-algebras.

\mypart{The iterated bar construction and iterated loop spaces}\label{CochainModels}

\renewcommand{\themainthm}{\arabic{mainthm}}
\setcounter{mainthm}{0}

In this concluding part,
we study applications of our results to cochain complexes of spaces
and iterated loop spaces.

To fix our framework,
a space $X$ refers to a simplicial set and we consider the normalized cochain complex $N^*(X)$
with coefficients in the ground ring $\kk$.
One proves that $N^*(X)$
can be equipped with the structure of a (unitary) $\EOp$-algebra, for some $E_\infty$-operad $\EOp$,
for all $X\in\Simp$,
so that the map $X\mapsto N^*(X)$
defines a functor from the category of simplicial sets $\Simp$
to the category of $\EOp$-algebras ${}_{\EOp}\C$
(see~\cite{HinichSchechtman} for a first proof of this result and~\cite{BergerFresse,McClureSmith} for more combinatorial constructions).
In the context of pointed spaces,
we replace $N^*(X)$ by the reduced cochain complex $\bar{N}^*(X)$
to use objects without unit.
Then we obtain that $\bar{N}^*(X)$
comes equipped with the structure of an $\EOp$-algebra, for some non-unitary $E_\infty$-operad $\EOp$,
in accordance with our conventions.

Let $F_X$ be any cofibrant replacement of $\bar{N}^*(X)$
in the model category of $\EOp$-algebras.
According to results of~\cite{Mandell},
the suspension $\Sigma F_X$ is equivalent to $\bar{N}^*(\Omega X)$
in the homotopy category of $\EOp$-algebras
provided that $\Omega X$ is connected and under standard finiteness and completeness assumptions on $X$
(see \cite[Theorem~1.2]{Mandell} and its proof in \emph{loc. cit.}).
For our needs,
we also have to record that the equivalence $\Sigma F_X\sim\bar{N}(\Omega X)$
is natural in the homotopy category of $\EOp$-algebras.

Theorem~\ref{BarStructure:thm:Existence} implies the existence of a well-defined iterated bar complex $B^n(A)$
for all $\EOp$-algebras $A$.
Theorem~\ref{HomotopyInterpretation:thm:HomotopyInterpretation}
implies that this iterated bar complex $B^n(A)$
is equivalent to the iterated suspension $\Sigma^n A$
if $A$ is a cofibrant $\EOp$-algebra.
Thus,
in the case of a cochain algebra $\bar{N}^*(X)$,
we obtain equivalences $B^n(\bar{N}^*(X))\sim B^n(F_X)\sim\Sigma^n F_X$
and we have $\Sigma^n F_X\sim\bar{N}^*(\Omega^n X)$
by an inductive application of the results of~\cite{Mandell}.
The assumptions which are made explicit in~\emph{loc. cit.} are reasonable for a single loop space $\Omega X$,
but give needlessly conditions in the case of higher iterated loop spaces,
at least in the context where the ground ring is a finite primary field $\kk = \FF_p$.
The actual purpose of this part is to review shortly the arguments of~\emph{loc. cit.}
and to examine assumptions on the space $X$
which ensure the equivalence $B^n(\bar{N}^*(X))\sim\bar{N}^*(\Omega^n X)$.

\medskip
One checks by a careful inspection of~\cite[\S 3, \S 5]{Mandell}
that the suspension $\Sigma F_X$ is equivalent to $\bar{N}^*(\Omega X)$ in the homotopy category of $\EOp$-algebras
as long as the cohomological Eilenberg-Moore spectral sequence of the path space fibration
\begin{equation*}
E_2 = \Tor^{H^*(X,\kk)}_*(\kk,\kk)\Rightarrow H^*(\Omega X,\kk)
\end{equation*}
converges.
By induction,
we obtain that the $n$-fold suspension $\Sigma^n F_X$
is equivalent in the homotopy category of $\EOp$-algebras
to $\bar{N}^*(\Omega^n X)$,
the cochain algebra of the $n$-fold iterated loop space of $X$,
if the cohomological Eilenberg-Moore spectral sequence of the path space fibration converges
for all loop spaces $\Omega^m X$, where $1\leq m\leq n$.
Record simply the next usual conditions
which ensure this convergence
in the context where the ground is either the rational field $\kk = \QQ$
or a finite primary field $\kk = \FF_p$:

\begin{mainfact}[See~\cite{Dwyer,Shipley}]\label{CochainModels:EilenbergMooreConvergence}
Let $n\geq 1$.
Suppose that:
\begin{enumerate}\renewcommand{\labelenumi}{\theenumi}
\renewcommand{\theenumi}{(1)${}_0$}
\item\label{EilenbergMooreConvergence:TrivialLowerHomotopyGroups}
The homotopy groups $\pi_*(X)$ are trivial for all $*\leq n$ (case $\kk = \QQ$).
\renewcommand{\theenumi}{(1)${}_p$}
\item\label{EilenbergMooreConvergence:FiniteLowerHomotopyGroups}
The homotopy groups $\pi_*(X)$ are finite $p$-groups for all $*\leq n$ (case $\kk = \FF_p$).
\end{enumerate}
and
\begin{enumerate}\renewcommand{\labelenumi}{\theenumi}
\renewcommand{\theenumi}{(2)}
\item\label{EilenbergMooreConvergence:CohomologyFiniteness}
The homotopy groups $\pi_*(X)$ are finitely generated in every degree $*>0$.
\end{enumerate}
Then the cohomological Eilenberg-Moore spectral sequence
\begin{equation*}
E_2 = \Tor^{H^*(\Omega^{m} X,\kk)}_*(\kk,\kk)\Rightarrow H^*(\Omega^{m+1} X,\kk)
\end{equation*}
converges for every $m<n$.
\end{mainfact}

The finiteness assumptions on homotopy groups imply that $H_*(\Omega^m X,\kk)$
forms a finitely generated $\kk$-module,
in every degree $*\geq 0$, for every $m\leq n$
(to check this folk claim, proceed by induction on the Postnikov tower of $\Omega^m X$).
Therefore the assumptions of Fact~\ref{CochainModels:EilenbergMooreConvergence}
imply the convergence of the Eilenberg-Moore spectral sequence by~\cite{Dwyer},
even if we deal with non-connected spaces by~\cite{Shipley}.

As a corollary, we obtain:

\begin{mainfact}\label{CochainModels:IteratedSuspension}
In the situations of Fact~\ref{CochainModels:EilenbergMooreConvergence},
we have a natural equivalence $\Sigma^n F_X\sim\bar{N}(\Omega^n X)$
in the homotopy category of $\EOp$-algebras.
\end{mainfact}

This assertion can also be proved by using models of Postnikov towers
in the category of $\EOp$-algebras (we also refer to~\cite{Mandell} for the definition of this model).
This finer argument would show that assumptions on lower homotopy groups $\pi_m(X)$, for $m<n$,
are unnecessary and can be dropped.
Thus the assertion of Fact~\ref{CochainModels:IteratedSuspension}
holds under the finiteness assumption~\ref{EilenbergMooreConvergence:CohomologyFiniteness}
of Fact~\ref{CochainModels:EilenbergMooreConvergence}
as long as the group $\pi_n(X)$ is trivial in the case $\kk = \QQ$,
a finite $p$-group in the case of a finite field $\kk = \FF_p$.

\medskip
As regards the iterated bar complex,
the existence of weak-equivalences $B^n(\bar{N}^*(X))\sim B^n(F_X)\sim\Sigma^n F_X$
and Fact~\ref{CochainModels:IteratedSuspension} imply:

\begin{mainthm}\label{thm:IteratedBarCompleteSpaceCase}
Under the assumptions of Fact~\ref{CochainModels:EilenbergMooreConvergence},
we have a natural isomorphism
\begin{equation*}
H^*(B^n\bar{N}^*(X))\simeq\bar{H}^*(\Omega^n X,\kk),
\end{equation*}
for every $n\geq 1$.\qed
\end{mainthm}

In the case $\kk = \FF_p$,
one can use the classical Bousfield-Kan tower $\{R_s X\}$
to improve Theorem~\ref{thm:IteratedBarCompleteSpaceCase}.

Recall simply that
$\bar{H}^*(X,\FF_p)\simeq\colim_s\bar{H}^*(R_s X,\FF_p)$ (see~\cite[Proposition III.6.5]{BousfieldKan} and \cite{Dror}).
Equivalently,
the natural morphism $\colim_s\bar{N}^*(R_s X)\rightarrow\bar{N}^*(X)$
defines a weak-equivalence in the category of $\EOp$-algebras.
Note that the bar complex commutes with sequential colimits
so that the natural morphism $\colim_s B^n(\bar{N}^*(R_s X))\rightarrow B^n(\bar{N}^*(X))$
defines a weak-equivalence as well.
Recall also that the spaces $R_s X$
satisfy assumptions \ref{EilenbergMooreConvergence:FiniteLowerHomotopyGroups} and~\ref{EilenbergMooreConvergence:CohomologyFiniteness}
of Fact~\ref{CochainModels:EilenbergMooreConvergence}
if the cohomology modules $H^*(X,\FF_p)$ are degreewise finite
(this folk assertion follows from a standard application of the spectral sequence of~\cite[\S X.6]{BousfieldKan}).

Hence Theorem~\ref{thm:IteratedBarCompleteSpaceCase} implies:

\begin{mainthm}\label{thm:IteratedBarTopologicalInterpretation}
Let $X$ be a pointed space whose cohomology modules $H^*(X,\FF_p)$ are degreewise finite.
Let $R_s X$ denote Bousfield-Kan' tower of $X$ (where $R = \FF_p$).
Then we have a natural isomorphism
\begin{equation*}
H^*(B^n\bar{N}^*(X))\simeq\colim_s\bar{H}^*(\Omega^n R_s X,\FF_p),
\end{equation*}
for every $n\geq 1$.\qed
\end{mainthm}

This result can be improved in good cases.
For instance,
if $X$ is a nilpotent space whose homotopy groups are degreewise finitely generated,
then theorems of~\cite{Shipley}
imply:
\begin{equation*}
\colim_s\bar{H}^*(\Omega^n_0 R_s X,\FF_p) = \bar{H}^*(\Omega^n_0 R_\infty X,\FF_p),
\end{equation*}
where $R_\infty X$ refers to Bousfield-Kan' $p$-completion of $X$
and $\Omega^n_0 Y$ denotes the connected component of the base point of $\Omega^n Y$,
for any pointed space $Y$.
Observe also that
$\Omega^n Y\sim\pi_n(Y)\times\Omega^n_0 Y$
and note that
\begin{equation*}
\colim_s\bar{H}^0(\Omega^n R_s X,\FF_p) = \colim_s \FF_p^{\pi_n(R_s X)} = \FF_p^{{\pi_n(R_\infty X)}^{\wedge}_p},
\end{equation*}
where the notation $\FF_p^{{\pi_n(R_\infty X)}^{\wedge}_p}$
refers to the module of maps $\alpha: \pi_n(R_\infty X)\rightarrow\FF_p$
which are continuous with respect to the $p$-profinite topology
(see~\cite[\S\S III-VI]{BousfieldKan}, see also~\cite{Morel} for a conceptual setting to do $p$-profinite topology).

\mypart{Acknowledgements}
I am grateful to several readers for many comments on the exposition of this work,
at various stages of preparation of the article.
Muriel Livernet and Haynes Miller motivated me to give a comprehensive survey of the background
in order to make the paper more accessible
and to improve the presentation of my results.
I thank Geoffrey Powell for a careful reading of a preliminary version of the article
and valuable observations on the writing.
I thank David Chataur, Mark Hovey, Jean-Louis Loday, Daniel Tanr\'e, for checking the final manuscript
and further accurate observations
which helped me to put the finishing touches.

\end{document}